\input{style/arxiv-general.cfg}
\documentclass[MSNbibl,number,citesort,rotating,dvips]{arxbj}
\makeatletter
   \@ifpackageloaded{graphicx}{}{\usepackage{graphicx}}
\makeatother
\usepackage{upgreek}
\usepackage{graphicx}

\aid{0}
\volume{21}
\issue{2}
\pubyear{2015}
\firstpage{1134}
\lastpage{1165}
\doi{10.3150/14-BEJ600} 

\makeatletter

\newcommand{\rright}{\right}
\newcommand{\lleft}{\left}
\newcommand{\rrvert}{\vert}
\newcommand{\llvert}{\vert}

\def\Prob{\operatorname{Prob}}
\def\Opt{\operatorname{Opt}}
\def\cN{{\mathcal{N}}}
\def\bR{{\mathbf{R}}}
\def\Z{{\mathbf{Z}}}
\def\e{\mathrm{e}}
\def\C{{\mathbf{C}}}
\def\R{{\mathbf{R}}}
\def\Prob{\operatorname{Prob}}
\def\Trace{\operatorname{Trace}}
\def\Risk{\operatorname{Risk}}
\def\cC{\mathcal{C}}
\def\cU{\mathcal{U}}
\newtheorem{theorem}{Theorem}[section]
\newtheorem{lemma}{Lemma}[section]
\newtheorem{proposition}{Proposition}[section]
\newremark{remark}{Remark}[section]

\def\cS{\mathcal{S}}
\def\bZ{{\mathbf{Z}}}
\def\mod{\operatorname{mod}}

\renewcommand{\emptyset}{\varnothing}

\newcommand{\Ceil}{\operatorname{Ceil}}
\newcommand{\Floor}{\operatorname{Floor}}

\def\sfrac#1#2{#1/#2}

\makeatother

\begin{document}
\begin{frontmatter}

\title{On detecting harmonic oscillations}
\runtitle{On detecting harmonic oscillations}

\begin{aug}
\author[A]{\inits{A.}\fnms{Anatoli} \snm{Juditsky}\corref{}\thanksref{A}\ead[label=e1]{juditsky@imag.fr}} \and
\author[B]{\inits{A.}\fnms{Arkadi} \snm{Nemirovski}\thanksref{B}\ead[label=e2]{nemirovs@isye.gatech.edu}}
\address[A]{LJK,
Universit\'e J. Fourier, B.P. 53, 38041 Grenoble
Cedex 9, France. \printead{e1}}
\address[B]{ISYE, Georgia Institute
of Technology, Atlanta, GA
30332, USA.\\ \printead{e2}}
\end{aug}

\received{\smonth{1} \syear{2013}}
\revised{\smonth{9} \syear{2013}}

%
\begin{abstract}
In this paper, we focus on the following testing problem: assume that
we are given observations of a
real-valued signal along the grid $0,1,\ldots,N-1$, corrupted by {white}
Gaussian noise. We want to
distinguish between two hypotheses: (a) the signal is a \textsl
{nuisance} -- a linear combination of
$d_n$ harmonic oscillations of {known} frequencies, and (b) signal is
the sum of a nuisance and a
linear combination of a given number $d_s$ of harmonic oscillations
with \textsl{unknown} frequencies,
and such that the distance (measured in the uniform norm on the grid)
between the signal and the set
of nuisances is at least $\rho>0$. We propose a computationally
efficient test for distinguishing
between (a) and (b) and show that its ``resolution'' (the smallest
value of $\rho$ for which (a) and
(b) are distinguished with a given confidence $1-\alpha$) is $\mathrm{O}
(\sqrt{\ln(N/\alpha)/N} )$,
with the hidden factor depending solely on $d_n$ and $d_s$ and
independent of the frequencies in question.
We show that this resolution, up to a factor which is polynomial in
$d_n, d_s$ and \textsl{logarithmic in $N$},
is the best possible under circumstances. We further extend the
outlined results to the case of nuisances
and signals \emph{close} to linear combinations of harmonic
oscillations, and provide illustrative numerical results.
\end{abstract}

%

\end{frontmatter}
%
\section{Introduction}
In this paper, we address the {following} detection problem. A signal
-- a two-sided sequence of reals $x=\{x_t,t=0,\pm1,\pm2,\ldots\}$ -- is
observed on the time horizon $0,1,\ldots,N-1$ according to
\[
y=x_0^{N-1}+\xi,
\]
where $\xi\sim\cN(0,I_N)$ is the white Gaussian noise and
$z_0^{N-1}=[z_0;\ldots ;z_{N-1}]$. Given $y$ we want to distinguish between
two hypotheses:
\begin{itemize}
\item\textsl{Nuisance hypothesis}: $x\in H_0$, where $H_0$ is
comprised of all \textsl{nuisances} -- linear combinations of $d_n$
harmonic oscillations of {known} frequencies.
\item\textsl{Signal hypothesis}: $x\in H_1(\rho)$, where $H_1(\rho
)$ is the set of all sequences $x$ representable as $s+u$ with the
``nuisance component'' $u$ belonging to $H_0$ and the ``signal
component'' $s$ being a sum of at most $d_s$ harmonic oscillations (of
whatever frequencies) such that the \emph{uniform} distance, on the
time horizon in question, from $x$ to all nuisance signals is at least
$\rho$:
\[
\min_{z\in H_0}\bigl\|x_0^{N-1}-z_0^{N-1}
\bigr\|_\infty\geq\rho.
\]
\end{itemize}
We are interested in a test which allows to distinguish, with a given
confidence $1-\alpha$, between the above two hypotheses for as small
``resolution'' $\rho$ as possible.

An approach to this problem which is generally advocated in the signal
processing literature is based on {the generalized likelihood ratio
test \cite{Han93,KH94,Dju96,NK11}. This test seems to enjoy the
optimal detection performance in the problem of distinguishing between
the ``pure noise'' hypothesis $H_0$ (no nuisance is present) and the
``signal hypothesis,'' $H_1$, that a signal which is a sum of $d_s$
sinusoids is present. For instance, under certain assumptions on the
signal frequencies, this test is claimed \cite{NK11} to distinguish
$1-\alpha$-reliably\footnote{We say that a test distinguishes
$1-\alpha$-reliably between the hypotheses, say, $H_0$ and $H_1$ if
its risk which is, in our case, the maximal probability of rejecting
the hypothesis when it is true, is bounded by~$\alpha$.} between $H_0$
and $H_1$ if the $\ell_2$-norm of the signal is $\geq c\sqrt{d_s\ln
(N/\alpha)}$. However, implementation of this detector requires
computing minimal log-likelihood -- the global optimal value of the
optimization problem
\[
\min_{a,\omega,\phi}\sum_{k=0}^{N-1}
\Biggl(y_k-\sum_{j=1}^{d_s}a_j
\sin(\omega_jk+\phi_j) \Biggr)^2
\]
(here the minimum is taken with respect to all problem parameters
$(a_j,\omega_j,\phi_j), j=1,\ldots,d_s$), and becomes numerically
challenging already for very moderate problem dimension $d_s$. To
circumvent numerical problems, under certain assumptions on the signal
frequencies (i.e., they are not ``too close to each other'') the test
can be applied ``sequentially'' \cite{NK11}, when ``fitting one
frequency at a time.'' A second family of tests, which does not require
estimation of unknown frequencies and amplitudes, relies upon} noise
subspace methods, such as multiple signal classification (MUSIC) \cite
{Pis73,mus1} (see also \cite{mus2,HanQui} for detailed presentation of
these techniques). To the best of our knowledge, no theoretical bounds
for the resolution of such tests are available. A~different test for
the case when no nuisance is present, based on the normalized
periodogram, has been proposed in \cite{Fish29}. The properties of
this test and of its various modifications were extensively studied in
the statistical literature (see, e.g., \cite
{Whittle,Han70,Chiu,HanQui,RB74,QK94}).
However, few theoretical results on the power of this test are
available for the case of sequence $x$ not being a linear combination
of Fourier harmonics $\mathrm{e}^{2\uppi\imath kt/N}$, $k=0,1,\ldots,N-1$ under
signal hypothesis. For instance, in the paper \cite{Davies87}, brought
to our attention by the referee, the properties of the periodogram test
are analysed in the problem of detection of one sinusoid, where it is
shown that asymptotically (when $N\to\infty$, the reliability of the
test is not ``too high,'' and the signal frequency is ``not too close''
to $0$ or $1$ but is otherwise arbitrary) the sensitivity of the test
can be improved by the factor up to ${\uppi\over2}$, if maximization of
the periodogram is carried over all frequencies in $[0,1]$ instead of\vspace*{1pt}
the set of ``Fourier frequencies'' $k\over N$, $k=0,\ldots,N-1$.\vadjust{\goodbreak}

In this paper, we show that a good solution to the outlined problem is
offered by an extremely simple test as follows.
\begin{quote}
Let $F_Nu= \{{1\over\sqrt{N}}\sum_{t=0}^{N-1}u_t\exp\{
2\uppi\imath kt/N\} \}_{k=0}^{N-1}\dvtx \C^N\to\C^N$ be the Discrete
Fourier Transform. Given the observation $y$, we solve the convex
optimization problem
\[
\Opt(y)=\min_z \bigl\{\bigl\|F_N
\bigl(y-z_0^{N-1} \bigr)\bigr\|_\infty\dvtx z\in
H_0 \bigr\}
\]
and compare the optimal value with a threshold $q_N(\alpha)$ which is
a valid upper bound on the $1-\alpha$-quantile of $\|F_N\xi\|_\infty
$, $\alpha\in(0,1)$ being a given tolerance:
\[
\Prob_{\xi\sim\cN(0,I_N)} \bigl\{\|F_N\xi\|_\infty>q_N(
\alpha ) \bigr\}\leq\alpha.
\]
If $\Opt(y)\leq q_N(\alpha)$, we accept the nuisance hypothesis,
otherwise we claim that a signal is present.
\end{quote}
It is immediately seen that the outlined test rejects the nuisance
hypothesis when it is true with probability at most $\alpha$.\footnote
{This fact is completely independent of what the nuisance hypothesis is
-- it remains true when $H_0$ is an arbitrary set in the space of
signals.} Our main result (Theorem~\ref{themain1}) states that \textsl
{the probability to reject the signal hypothesis when it is true is
$\leq\alpha$, provided that the resolution $\rho$ is not too small},
specifically, {for an appropriately chosen universal function $\cC
(\cdot)$,}
{\renewcommand{\theequation}{!}
\begin{equation}
\rho\geq\cC(d_n+d_s)\sqrt{\ln(N/\alpha)/N}.
\end{equation}}
{In the simplest case when the zero hypothesis just states that $x=0$,
the above test is similar to the standard tests based on the maximum of
the periodogram; when $H_0$ is nontrivial, our test can be seen as a
natural modification of the classical construction.}

Some comments are in order.
\begin{itemize}[B.]
\item[A.] {Our principal} contributions, as we see them, are in
\begin{itemize}[--]
\item[--] deriving theoretical upper (and close to them lower) bounds on
the resolution of the test under quite general assumptions on nuisances
and signals of interest; it should be stressed that these bounds depend
solely on the cardinalities of the sets of ``participating''
frequencies. Moreover, we allow for ``multiplicities of frequencies,''
so that our ``sum of harmonic oscillations'' can be an algebraic
polynomial or a sum of products of algebraic polynomials and harmonic
oscillations. We are not aware of comparable, in terms of generality,
existing theoretical results on the performance of spectral analysis tests;
\item[--] demonstrating that in order to achieve nearly optimal resolution,
one can restrict the periodogram to the frequencies participating in
the Discrete Fourier Transform (the traditional recommendation is to
consider denser frequency sets).
\end{itemize}
\item[B.] We show that the power of our test when applied to
the detection problem in question is nearly as good as it can be:
precisely, for every pair $d_n$, $d_s$ and properly selected $H_0$,
no test can distinguish $(1-\alpha)$-reliably between $H_0$ and
$H_1(\rho)$ when $\rho<\mathrm{O}(1)d_s\sqrt{\ln(1/\alpha)/N}$. Here and
from now on, $\mathrm{O}(1)$'s are appropriately chosen positive absolute constants.\vadjust{\goodbreak}
\item[C.] We are measuring the resolution in the ``weakest''
of all natural scales, namely, via the \textsl{uniform} distance from
the signal to the set of nuisances.
When passing from\vspace*{1pt} the uniform norm to the normalized Euclidean norm
$|x_0^{N-1}|_2:=\|x_0^{N-1}\|_2/\sqrt{N}\leq\|x_0^{N-1}\|_\infty$,
an immediate lower bound on the resolution
which allows for reliable detection becomes $\mathrm{O}(1)\sqrt{\ln(1/\alpha
)/N}$. In the case when, as in our setting, signals obeying $H_0$ and
$H_1(\rho)$ admit parametric description involving $K$ parameters,
this lower bound, up to a factor logarithmic in $N$ and linear in $K$,
is also an upper resolution bound, and the associated test is based on
estimating the Euclidean distance from the signal underlying the
observations to the nuisance set. Note that, in general, the $|\cdot
|_2$-norm can be smaller than $\|\cdot\|_\infty$ by a factor as large
as $\sqrt{N}$, and the fact that ``energy-based'' detection allows to
distinguish well between parametric hypotheses ``separated'' by
$\mathrm{O}(\sqrt{\ln(N/\alpha)/N})$ in $|\cdot|_2$ norm does \textsl{not}
automatically imply the possibility to distinguish between hypotheses
separated by $\mathrm{O}(\sqrt{\ln(N/\alpha)/N})$ in the uniform
norm.\footnote{Indeed, let $H_0$ state that the signal is 0, and
$H_1(\rho)$ state that the signal is $\geq\rho$ at $t=0$ and is zero
for all other $t$'s. These two hypotheses cannot be reliably
distinguished unless $\rho\geq \mathrm{O}(1)$, that is, the $\|\cdot\|_\infty
$ resolution in this case is much larger than $\mathrm{O}(\sqrt{\ln(N/\alpha
)/N})$.} The latter possibility exists in the situation we are
interested in due to the particular structure of the specific nuisance
and signal hypotheses;
this structure allows also for a dedicated \textsl{non-energy-based} test.
\item[D.] For the sake of definiteness, throughout the paper
we assume that the observation noise is the standard white Gaussian
one. This assumption is by no means critical: on a straightforward
inspection of what follows, \textsl{whatever be the observation noise,
with $q_N(\alpha)$ defined as (an upper bound on) the $(1-\alpha
)$-quantile of $\|F_N\xi\|_\infty$,
the above test $(1-\alpha)$-reliably distinguishes between the
hypotheses $H_0$ and $H_1(\rho)$, provided that $\rho\geq\cC
(d_n+d_s)q_N(\alpha)/\sqrt{N}$.} For example, the results of Theorems
\ref{themain1} and \ref{themain11} remain valid when the observation
noise is of the form $\xi=\{\xi_t=\sum_{\tau=-\infty}^\infty
\gamma_\tau\eta_{t-\tau}\}_{t=0}^{N-1}$ with deterministic $\gamma
_\tau$, $\sum_\tau|\gamma_\tau|\leq1$, and independent $\eta
_t\sim\cN(0,1)$.
\item[E.] The main observation underlying the results on the
resolution of the above test is as follows: \textsl{when $x$ is the
sum of at most $d$ harmonic oscillations, $\|F_Nx\|_\infty\geq\cC
(d)\sqrt{N}\|x_0^{N-1}\|_\infty$ with some universal positive
function $\cC(d)$.\footnote{Since $\|x\|_2=\|F_Nx\|_2\geq\|F_Nx\|
_\infty$, it follows that for the aforementioned $x$ one has $\cC
(d)^{-1}\|x\|_\infty\leq N^{-1/2}\|x\|_2\leq\|x\|_\infty$ -- an
important by its own right fact which we were unable to find in
mathematical literature.}} This observation originates from \cite
{Nem81} and, along with its modifications and extensions, was utilized,
for the time being in the denoising setting, in \cite
{Nem92,GoNem97,JuNem1,JuNem2}. It is worth mentioning that it also
allows to extend, albeit with degraded constants, the results of
Theorems \ref{themain1}
and \ref{themain11} to multi-dimensional setting.
\end{itemize}

The rest of this paper is organized as follows. In Section~\ref{sect:Problem}, we give a detailed description of the detection
problems $(P_1)$, $(P_2)$, $(N_1)$ and $(N_2)$, we are interested in
(where $(P_2)$ is the problem we have discussed so far).
Our test is presented in Section~\ref{sectTest} where we also provide
associated resolution bounds for these problems. Next, in Section~\ref{sectLower}, we present lower bounds on ``good'' (allowing for
$(1-\alpha)$-reliable hypotheses testing) resolutions, while in
Section~\ref{sectNumres} we describe some numerical illustrations. The
proofs of results of Sections~\ref{sectTest} and \ref{sectLower} are
put into Section~\ref{sec:proofs}.

\section{Problem description}\label{sect:Problem}
Let $\cS$ stand for the space of all two-sided real sequences $z=\{
z_t\in\bR\}_{t=-\infty}^\infty$. Assume that a discrete time signal
$x\in\cS$ is observed on the time horizon $0\leq t< N$ according to
\setcounter{equation}{0}
\begin{equation}
\label{eqobs} y=x_0^{N-1}+\xi, \qquad \xi\sim
\cN(0,I_N),
\end{equation}
where (and from now on) for $z\hspace*{-0.3pt}\in\hspace*{-0.3pt}\cS$ and integers $p\hspace*{-0.3pt}\leq\hspace*{-0.3pt} q$, $z_p^q$
stands for the vector $[z_p;z_{p+1};\ldots ;z_q]$.

In the sequel, we are interested in the case when the signal is a
linear combination of a given number
of harmonic oscillations. Specifically, let $\Delta$ stand for the
shift operator on $\cS$:
\[
(\Delta z)_t=z_{t-1},\qquad z\in\cS.
\]
Let $\Omega_d$ be the set of all unordered collections ${\mathbf
{w}}=\{\omega_1,\ldots,\omega_d\}$ of $d$ reals which are ``symmetric
$\mod 2
\uppi$,'' meaning that for every $a$, the number of $\omega_i$'s equal,
modulus $2\uppi$, to $a$ is exactly the same as the number of $\omega
_i$'s equal, modulus $2\uppi$, to $-a$. We associate with ${\mathbf
{w}}\in\Omega_d$ the real algebraic polynomial
\[
p_{\mathbf{w}}(\zeta)=\prod_{\ell=1}^d
\bigl(1-\exp\{\imath\omega_\ell \}\zeta \bigr)
\]
and the subspace $\cS[{\mathbf{w}}]$ of $ \cS$, comprised of $x\in
\cS$ satisfying the homogeneous finite-difference equation
%
\begin{equation}
\label{equation} p_{\mathbf{w}}(\Delta)x\equiv0.
\end{equation}
In other words, $\cS[\{\omega_1,\ldots,\omega_d\}]$ is comprised of all
real two-sided sequences of the form
\[
x_t=\sum_{\ell=1}^d
\bigl[p_\ell(t)\cos(\omega_\ell t)+q_\ell (t)\sin(
\omega_\ell t) \bigr]
\]
with real algebraic polynomials $p_\ell(\cdot)$, $q_\ell(\cdot)$ of
degree $<m_\ell$, where $m_\ell$ is the multiplicity, $\mod 2\uppi$,
of $\omega_\ell$ in $\mathbf{w}$. We set
\[
\cS_d=\bigcup_{{\mathbf{w}}\in\Omega_d}\cS[{\mathbf{w}}].
\]

\begin{remark} \label{rem1} In what follows, we refer to the
reals $\omega_i$ constituting ${\mathbf{w}}\in\Omega_d$ as the
\textsl{frequencies} of a signal from $\cS[{\mathbf{w}}]$. A reader
would keep in mind that the number of ``actual frequencies'' in such a
signal can be less than $d$: for instance, frequencies in ${\mathbf
{w}}$ different from $0 \mod 2\uppi$ and $\uppi\mod 2\uppi$ go in
``symmetric pairs'' $(\omega,\omega'=-\omega\mod  2\uppi)$, and such a
pair gives rise to a single ``actual frequency.''
\end{remark}

Given a positive integer $N$, real $\epsilon\geq0$, and ${\mathbf
{w}}\in\Omega_d$, we set
\[
\cS^{N,\epsilon}[{\mathbf{w}}]= \bigl\{x\in\cS\dvtx \bigl\| \bigl[p_{\mathbf
{w}}(
\Delta)x \bigr]_0^{N-1}\bigr\|_\infty\leq\epsilon \bigr
\}.
\]
Finally, we denote by $\cS^{N,\epsilon}_d$ the set
\[
\cS^{N,\epsilon}_d=\bigcup_{{\mathbf{w}}\in\Omega_d}\cS
^{N,\epsilon}[{\mathbf{w}}].
\]
When $N$ is clear from the context, we shorten the notations $\cS
^{N,\epsilon}[{\mathbf{w}}]$, $\cS^{N,\epsilon}_d$ to $\cS
^{\epsilon}[{\mathbf{w}}]$ and $\cS^{\epsilon}_d$, respectively.

In the definitions above, it was tacitly assumed that $d$ is a positive
integer. It makes sense to allow also for the case of $d=0$. By
definition, $\Omega_0$ is comprised of the empty collection ${\mathbf
{w}}=\emptyset$ and $p_\emptyset(\zeta)\equiv1$. With this
convention, $\cS^{N,\epsilon}[\emptyset]=\{x\in\cS\dvtx \|x_0^{N-1}\|
_\infty\leq\epsilon\}$.

Observe that the family $\cS^{N,\epsilon}_d$, $d\ge2$, is quite
rich. For instance, it contains ``smoothly varying signals'' (case of
$w_i=0 \mod 2\uppi$), along with ``fast varying'' -- amplitude-modulated
and frequency-modulated signals (see \cite{GoNem97,JuNem1} for more examples).

We detail now the hypothesis testing problems about the sequence $x$
via observation $y$ given by (\ref{eqobs}).
In what follows $d_s$ and $d_n$ are given positive integers, and $\rho
$, $\epsilon_n$, $\epsilon_s$ are given positive reals.
\begin{enumerate}[$(N_2)$]
\item[$(P_1)$] The ``basic'' hypothesis testing problem we consider is
that of testing of a simple nuisance hypothesis $\{x=0\}$ against the
alternative that a signal $x\in\cS_{d_s}$ ``is present,'' meaning
that the uniform norm of the signal on the observation window
$[0,\ldots,N-1]$ exceeds certain threshold $\rho>0$. In other words, we
consider the following set of hypotheses:
\begin{eqnarray*}
H_0&=&\{x= 0\},
\\
H_1(\rho)&=& \bigl\{x\in\cS_{d_s}\dvtx
\bigl\|x_0^{N-1} \bigr\|_\infty\geq\rho \bigr\}.
\end{eqnarray*} %
\item[$(P_2)$] We suppose that $x\in\cS$ decomposes into ``signal''
and ``nuisance'':
\[
x=s+u,
\]
where $s$ is the signal of interest and a \textsl{nuisance} $u$
belongs to a subspace $\cS[\overline{\mathbf{w}}]$, assumed to be
known a priori. We consider a composite nuisance hypothesis
that $x$ is a ``pure nuisance,'' and the alternative (signal
hypothesis) that useful signal $s$ does not vanish, and the deviation,
when measured in the uniform norm on the observation window, of
``signal$\,+\,$nuisance'' from the nuisance subspace is at least $\rho>0$.
Thus, we arrive at the testing problem: given $\overline{\mathbf
{w}}\in\Omega_{d_n}$ decide between the hypotheses
\begin{eqnarray*}
H_0&=& \bigl\{x=u\in\cS[\overline{\mathbf{w}}] \bigr\},
\\
H_1(\rho)&=& \Bigl\{ %
x=u+s\dvtx u\in \cS[\overline{
\mathbf{w}}], s\in\cS_{d_s},
\\
&&\hphantom{ \bigl\{}\mbox{such that } \min_{z} \bigl\{
\bigl\|[x-z]_0^{N-1} \bigr\|_\infty\dvtx z\in \cS[\overline{
\mathbf{w}}] \bigr\}\geq\rho \Bigr\}.
\end{eqnarray*}
Clearly, problem $(P_1)$ is a particular case of $(P_2)$ with $d_n=0$
(and thus $\cS_{d_n}=\{0\}$ is a singleton).

\item[$(N_1)$] Given $\epsilon_n>0$ and $\overline{\mathbf{w}}\in
\Omega_{d_n}$, consider the nonparametric nuisance hypothesis that the
nuisance $u\in\cS^{N,\epsilon_n}[\overline{\mathbf{w}}]$ with some known
$\overline{\mathbf{w}}$. The signal hypothesis is that the useful
signal $s\in\cS_{d_s}$ is present, and $x=s+u$ deviates from the
nuisance set on the observation window by at least $\rho>0$ in the
uniform norm:
\begin{eqnarray*}
H_0&=& \bigl\{x=u\in\cS^{N,\epsilon_n}[\overline{\mathbf{w}}] \bigr
\},
\\
H_1(\rho)&=& \Bigl\{ %
x=u+s\dvtx u\in \cS^{N,\epsilon_n}[
\overline{\mathbf{w}}], s\in\cS _{d_s},
\\
&&\hphantom{ \bigl\{}\mbox{such that } \min_{z} \bigl\{
\bigl\|[x-z]_0^{N-1} \bigr\|_\infty\dvtx z\in
\cS^{N,\epsilon_n}[\overline{\mathbf{w}}] \bigr\}\geq\rho \Bigr\}.
\end{eqnarray*}
\item[$(N_2)$] The last decision problem is a natural extension of
$(N_1)$: we consider the problem of testing a nonparametric nuisance
hypothesis against a nonparametric signal alternative that the useful
signal $s\in\cS^{N,\epsilon_s}_{d_s}$ is present:
\begin{eqnarray*}
H_0&=& \bigl\{x=u\in\cS^{N,\epsilon_n}[\overline{\mathbf{w}}] \bigr
\},
\\
H_1(\rho)&=& \Bigl\{ %
x=u+s\dvtx u\in \cS^{N,\epsilon_n}[
\overline{\mathbf{w}}], s\in\cS ^{N,\epsilon_s}_{d_s},
\\
&&\hphantom{ \bigl\{}\mbox{and such that } \min_{z} \bigl\{
\bigl\|[x-z]_0^{N-1} \bigr\|_\infty\dvtx z\in
\cS^{N,\epsilon_n}[\overline{\mathbf{w}}] \bigr\}\geq\rho \Bigr\}.
\end{eqnarray*}
Note that problem ($N_1$) is a particular case of ($N_2$) with
$\epsilon_s=0$.
\end{enumerate}
In the sequel,
we refer to the sequences obeying $H_0$ (resp., $H_1=H_1(\rho)$) as
\textsl{nuisance} (resp., \textsl{signal}) sequences.

Let $\varphi(\cdot)$ be a test, that is, a Borel function on $\bR^N$
taking values in $\{0,1\}$, which receives on input observation (\ref
{eqobs}) (along with {parameters describing} $H_0$ and $H_1$).
The event\linebreak[4]  $\{\varphi(y)=1\}$ corresponds to
rejecting
the hypothesis $H_0$, while $\{\varphi(y)=0\}$ implies that $H_1$
is rejected. The quality of the test
is characterized by the error probabilities -- the probabilities of
rejecting erroneously each of the hypotheses:
\[
\varepsilon_0(\varphi; H_0)=\sup_{x\in H_0}
\Prob_x \bigl\{\varphi(y)= 1 \bigr\}, \qquad \varepsilon_1
\bigl( \varphi; H_{1}(\rho) \bigr)= \sup_{x\in H_{1}(\rho)}
\Prob_x \bigl\{\varphi(y)= 0 \bigr\}.
\]
We define the \emph{risk of the test} as
\[
\Risk(\varphi,\rho)=\max \bigl\{ \varepsilon_0(\varphi;
H_0), \varepsilon_1 \bigl(\varphi; H_{1}(
\rho ) \bigr) \bigr\}.
\]
Let $\alpha\in(0,1/2)$ be given. In this paper, we address the
following question: \textsl{for the testing problems above, what is
the smallest possible $\rho$ such that one can distinguish $(1-\alpha
)$-reliably between the hypotheses $H_0$ and $H_1=H_1(\rho)$ via
observation} (\ref{eqobs}) (i.e., is such that $\Risk(\varphi,\rho)\leq\alpha$).
In the sequel, we refer to such $\rho$ as to \emph{$\alpha
$-resolution} in the testing problem in question, and our goal is to
find reasonably tight upper and lower bounds on this resolution along
with the test underlying the upper bound.\vspace*{-1pt}
%
\section{Basic test and upper resolution bounds}
\label{sectTest}
In this section, we present a simple test which provides some upper
bounds on the $\alpha$-resolutions in problems ($P_1$)--($N_2$).

Let $\Gamma_N=\{\mu_\tau=\exp\{2\uppi\imath\tau/N\}\dvtx 0\leq\tau
<N-1\}$ be the set of all roots $\mu\in\C$ of unity of degree $N$,
and let $F_N\dvtx \C^N\to\C(\Gamma_N)$ be the normalized Fourier transform:
%
\begin{eqnarray}
\label{eqFour} [F_Nf](\mu)={1\over\sqrt{N}}\sum
_{t=0}^{N-1}f_t\mu^t,\qquad
\mu\in \Gamma_N.
\end{eqnarray}
Note that (\ref{eqFour}) can also be seen as a mapping from $\cS$ to
$\C
(\Gamma_N)$.

Given a tolerance $\alpha\in(0,1/2)$, let $q_N(\alpha)$ be 
the $(1-\alpha)$-quantile of $\|F_N\xi\|_\infty$, so that
\[
\Prob_{\xi\sim\cN(0,I_N)} \bigl\{\|F_N\xi\|_\infty\geq
q_N(\alpha) \bigr\} \leq\alpha.
\]
Let $Q_{\cN}(\alpha)$ be the $(1-\alpha)$-quantile of the standard
normal distribution:
\[
{1\over\sqrt{2\uppi}}\int_{Q_{\cN}(\alpha
)}^\infty \mathrm{e}^{-{s^2/2}}\,\mathrm{d}s=\alpha.
\] In the sequel, we use the
following immediate bound for $q_N(\cdot)$:\footnote{{For instance,
for even $N$, $[F_N\xi](1)$ and $[F_N\xi](-1)$ are (real) standard
normal variables, and $|[F_N\xi](\mu)|^2$ follows exponential
distribution when $\mu\neq\pm1$. Indeed, in the latter case the
joint distribution of real and imaginary parts of $[F_N\xi](\mu)$ is
$\cN(0, \frac{1}{2}I_2)$, so that $2[F_N\xi](\mu)|^2\sim\chi_2^2$
($\chi^2$ distribution with 2 degrees of freedom). As a result, $\|
F_N\xi
\|_\infty$ can be bounded with $\max (|F_N\xi](1)|, |F_N\xi
](-1)|,\sqrt{\eta_{N-2}} )$, where $\eta_{N-2}$ is the maximum
of $N-2$ independent exponential random variables, what implies the
bound (\ref{onecantake}) in this case.

It is worth to mention that the same argument leads to the bound for
$q_N(\alpha)$ which is equivalent to (\ref{onecantake})\vspace*{1pt} in the case when
the observation noise is of the form $\xi=\{\xi_t=\sum_{\tau
=-\infty}^\infty\gamma_\tau\eta_{t-\tau}\}_{t=0}^{N-1}$ with
deterministic $\gamma_\tau$, $\sum_\tau|\gamma_\tau|\leq1$, and
independent $\eta_t\sim\cN(0,1)$. Indeed, in this case we have
$[F_N\xi](\mu)=\sum_{\tau=-\infty}^{\infty} \overline{\gamma
}_{\mu,\tau} \eta_\tau$, with $\sum_{\tau} |\overline{\gamma
}_{\mu,\tau} |^2\leq1$. In other words,
the joint distribution of real and imaginary parts of $[F_N\xi](\mu)$
is $\cN(0,\Sigma_2)$, with $\Trace(\Sigma_2)\le1$. Then for any
$\rho>0$,
\[
\Prob \bigl\{\|F_N\xi\|_\infty\ge\rho \bigr\}\le N\max
_{\mu\in
\Gamma_N}\Prob \bigl\{\bigl|[F_N\xi](\mu)\bigr|\ge\rho \bigr\}
\le2N\Prob \bigl\{|\zeta|\ge\rho \bigr\},
\]
where $\zeta\sim\cN(0,1)$, and $\Prob \{\|F_N\xi\|_\infty\ge
Q_\cN ({\alpha\over4N} ) \}\le\alpha$.
}
}
%
\begin{eqnarray}
\label{onecantake} q_N(\alpha)&\le&\lleft\{ %
\begin{array} {l@{
\qquad}l} \displaystyle \inf_{0\le s\le1} \max \biggl[{Q_{\cN}} \biggl(
\displaystyle {s\alpha\over4} \biggr), \sqrt{\ln \biggl(\displaystyle {N-2\over2(1-s)\alpha} \biggr)}
\biggr]&\mbox{for $N$ even},
\nonumber
\\
\displaystyle \inf_{0\le s\le1} \max \biggl[{Q_{\cN}} \biggl(
\displaystyle {s\alpha\over
2} \biggr), \sqrt{\ln \biggl(\displaystyle {N-1\over2(1-s)\alpha} \biggr)}
\biggr]&\mbox{for $N$ odd} \end{array} %
\rright.
\nonumber
\\[-8pt]
\\[-8pt]
&\asymp&\sqrt{\ln ({N/ \alpha} )},
\nonumber
\end{eqnarray}
where $a\asymp b$ means that the ratio $a/b$ is in-between absolute
positive constants.

\textit{The test} we are about to consider (and which we refer to
as \emph{basic test} in the sequel) is as follows:
\begin{enumerate}[2.]
\item[1.] Given $y$, we solve the convex optimization problem
%
\begin{equation}
\label{convexproblem0} \Opt_\mathcal{Z}(y)=\min_{z\in\mathcal{Z}}
\bigl\|F_N \bigl(y-z_0^{N-1} \bigr)\bigr\| _\infty,
\end{equation}
where the set $\mathcal{Z}$ is defined according to
%
\begin{eqnarray}
\mathcal{Z}=\lleft\{ %
\begin{array} {l@{\qquad}l} \{0\}&\mbox{for
problem $(P_1)$},
\\
\cS[\overline{\mathbf{w}}]&\mbox{for problem $(P_2)$},
\\
\cS^\epsilon[\overline{\mathbf{w}}]&\mbox{for problems
($N_1$) and ($N_2$)}. \end{array} %
\rright.
\label{convexproblem1}
\end{eqnarray}
\item[2.] We compare $\Opt_\mathcal{Z}(y)$ to $q_N(\alpha)$, where
$\alpha$ is a given tolerance: if ${\Opt_\mathcal{Z}(y)}\leq
q_N(\alpha)$, we accept $H_0$, otherwise we accept $H_1$.
\end{enumerate}
We describe now the properties of the basic test as applied to problems
$(P_1)$, $(P_2)$, $(N_1)$ and $(N_2)$.

\begin{theorem}\label{themain1} The risk of the basic test as applied
to problems $(P_1)$, $(P_2)$ is bounded by $\alpha$, provided that
$d_*=d_n+d_s>0$ and
%
\begin{equation}
\label{resolution} \rho\geq \mathrm{O}(1)d_*^3\ln(2d_*) q_N(
\alpha)N^{-1/2}=\mathrm{O}(1) d_*^3\ln(2d_*)
\sqrt{N^{-1}\ln ({N/\alpha} )}
\end{equation}
with properly chosen positive absolute constants $\mathrm{O}(1)$.\footnote
{{Note that the signal $s\in\cS_{d_s}$ is assumed to be at a distance
$\rho>0$ from the nuisance subspace (and thus from the origin). In
other words, under the premise of the Theorem~\ref{themain1} $d_s>0$.}}
\end{theorem}

The result for the nonparametric problems ($N_1$) and ($N_2$) is similar.

%
\begin{theorem}\label{themain11} The risk of the basic test as applied
to problems $(N_1)$, $(N_2)$ is bounded by $\alpha$, provided that
$d_*=d_n+d_s>0$,
$\rho$ satisfies (\ref{resolution}) with properly selected $\mathrm{O}(1)$\textup{'}s
and, in addition, $\epsilon_n$ and $\epsilon_s$ are small enough,
specifically,
%
\begin{equation}
\label{epsilonsaresmall} N^{d_n+\sfrac{1}{2}}\epsilon_n+N^{d_s+\sfrac{1}{2}}
\epsilon_s\leq \mathrm{O}(1)q_N(\alpha)
\end{equation}
with properly selected positive absolute constant $\mathrm{O}(1)$.
\end{theorem}

The proofs of Theorems \ref{themain1} and \ref{themain11} are
relegated to Section~\ref{sec:proofs}.

Theorems \ref{themain1} and \ref{themain11} provide us with upper
resolution bounds independent of the frequencies constituting
$\overline{\mathbf{w}}$ and ${\mathbf{w}}$.
When $\epsilon_n$, $\epsilon_s$ are ``small enough,'' so that (\ref
{epsilonsaresmall}) holds true (we refer to the corresponding range of
problems' parameters as the \emph{parametric zone}), our upper bound
on $\alpha$-resolution in all testing problems of interest is
essentially the same as in the case of $\epsilon_n=\epsilon_s=0$~--
it is $\cC(d_n+d_s)\sqrt{\ln(N/\alpha)/N}$ with the factor $\cC
(d)=\mathrm{O}(1)d^3\ln(2d)$ depending solely on $d$.

On the other hand, when $\epsilon_n$ and $\epsilon_s$ are \textsl
{not} ``small enough,'' that is, when $(N,\epsilon_n,\epsilon_s)$
are not in the range described by (\ref{epsilonsaresmall}), and $\rho
=\mathrm{O}(1)\sqrt{\ln(1/\alpha)/N}$, $(1-\alpha)$-reliable decision
between the hypotheses participating in $(N_1)$, $(N_2)$ becomes
impossible, whatever be the test (see items (ii) in Propositions \ref
{lowbndC}, \ref{lowbndD} below). {This phenomenon is by no means
surprising}: indeed, allowing for $\epsilon_n$ and/or $\epsilon_s$ to
be positive means allowing for \textsl{nonparametric} components in
the nuisance and/or nonnuisance signals.\footnote{Think, for example,
of a ``{plain} nonparametric'' signal {$x_0^{N-1}$, which is a
restriction on the grid $\{i/N\}_{i=0}^{N-1}$ of a smooth function on
$[0,1]$. It is immediately seen that such a signal} belongs to $\cS
^{N,\epsilon}[0]$ with $\epsilon\le\|f'\|_\infty/N$.} When these
components are not ``small enough,'' their presence changes
significantly the resolution level at which reliable test is possible.
The study of problems $(N_1)$ and $(N_2)$ in the nonparametric range
(i.e., with $(N,\epsilon_n,\epsilon_s)$ beyond the range given by
(\ref{epsilonsaresmall})) goes beyond the scope of this paper\footnote
{For ``de-noising'' analogies of our testing problems in the
nonparametric range, see \cite{Nem92,GoNem97,JuNem1,JuNem2}.} and
deserves, we believe, a dedicated study. Our ``educated guess'' is that
correct nonparametric version of the basic test in problems $(N_1)$,
$(N_2)$ should be as follows: given $(N,\epsilon_n,\epsilon_s)$
\textsl{not} satisfying (\ref{epsilonsaresmall}), find the largest
integer $\overline{N}$ such that $(\overline{N},\epsilon_n,\epsilon
_s)$ does satisfy (\ref{epsilonsaresmall}), and then run the basic
test on, say, all $\overline{N}$-point subintervals of the $N$-point
observation horizon, inferring the validity of the nuisance hypothesis
if and only if it was accepted by the basic test on every
subinterval.\footnote{With this approach, the basic test on a
subinterval should be associated with confidence $\mathrm{O}(\alpha/N)$ instead
of $\alpha$, in order to account for the possibility of ``large
deviations'' in at least one of $N-\overline{N}+1$ runs of the basic test.}

\section{Lower resolution bounds}\label{sectLower}
The lower resolution bounds of this section complement the upper bounds
of Section~\ref{sectTest}.
We start with the parametric setting $(P_1)$ and $(P_2)$. Through this
section, $c_i(d_n,d_s)$ are properly selected positive and monotone
functions of their arguments.

\begin{proposition}\label{lowbndAB} Given integers $d_n\ge0$,
$d_s\geq1$, and a real $\alpha\in(0,1/2)$, consider problems $(P_1)$
and $(P_2)$ with {parameters} $d_n$ ($d_n=0$ in the case of problem
$(P_1)$), $d_s$, $\alpha$ and
$\overline{\mathbf{w}}=\{\overbrace{0,\ldots,0}^{d_n}\}$. Then for
properly selected $c_0(d_n,d_s)$ and all $N\geq c_0(d_n,d_s)$
the $\alpha$-resolution $\rho$ in the problems $(P_1)$ and $(P_2)$
admits the lower bound
\[
\mathrm{O}(1) d_s\sqrt{\ln(1/\alpha)/N}.
\]
\end{proposition}

We see that in the problem $(P_1)$ $\alpha$-resolution grows with
$d_s$ at least linearly. Note that by Theorem~\ref{themain1}, this
growth is at most cubic (more precisely, it is not faster than
$\mathrm{O}(1)d_s^3\ln(d_s)$). Besides this, we see that the upper bounds on
$\alpha$-resolution for problems $(P_1)$ and $(P_2)$ stemming from
Theorem~\ref{themain1} and associated with the basic test coincide,
within a factor depending solely on $d_n,d_s,N$ and logarithmic in $N$,
with lower bounds on $\alpha$-resolution.

We have the following lower bound on the $\alpha$-resolution in the
problem $(N_1)$.

\begin{proposition}\label{lowbndC}
Given integers $d_n>0$, $d_s\geq2$ and reals $\alpha\in(0,1/2)$,
${\epsilon_n}\geq0$, consider problem $(N_1)$ with {parameters}
$d_n$, $d_s$, $N$, ${\epsilon_n}$, $\alpha$ and
$\overline{\mathbf{w}}=\{\overbrace{0,\ldots,0}^{d_n}\}$.
Then for properly selected $c_i(d_n,d_s)>0$ depending solely on
$d_n,d_s$ and for all $N$ satisfying
%
\begin{eqnarray}
\label{proplet} N\ge c_0(d_n,d_s),
\end{eqnarray}
the $\alpha$-resolution $\rho_*(\alpha)$ in the problem $(N_1)$ satisfies:
\begin{enumerate}[(ii)]
\item[(i)] in the range $0\le{\epsilon_n}\le
c_1(d_n,d_s)N^{-d_n-1/2}\sqrt{\ln(1/\alpha)}$,
\begin{eqnarray*}
\rho_*(\alpha)\geq c_2(d_n,d_s)\sqrt{
\ln(1/\alpha)/N};
\end{eqnarray*}
\item[(ii)] in the range
%
\begin{eqnarray}
\label{range}&& c_3(d_n,d_s)N^{-{d_n}-1/2}
\sqrt{\ln(1/\alpha)}\leq{\epsilon_n}\leq c_1(d_n,d_s)N
\sqrt{\ln(1/\alpha)},
\\
&&\rho_*(\alpha) \geq c_4(d_n,d_s)
\bigl[{ \epsilon_n} N^{{d_n}+1/2} \bigl[\ln(1/\alpha)
\bigr]^{-1/2} \bigr]^{{1/(2{d_n}+3)}}\sqrt {N^{-1}\ln(1/\alpha)}.
\nonumber
\end{eqnarray}
\end{enumerate}
\end{proposition}

In the case of the problem $(N_2)$, we have a similar lower bound on
$\alpha$-resolution when $\epsilon_n\le\epsilon_s$.

\begin{proposition}\label{lowbndD}
Given integers $d_n>0$, $d_s\geq2$ and reals $\alpha\in(0,1/2)$,
${\epsilon_n}\geq0$, consider problem $(N_2)$ with {parameters}
$d_n$, $d_s$, $N$, ${\epsilon_n}$, $\epsilon_s$, $\alpha$ and
$\overline{\mathbf{w}}=\{\overbrace{0,\ldots,0}^{d_n}\}$.

Assume that $0\le\epsilon_n\le\epsilon_s$.
Then for properly selected $c_i(d_n,d_s)>0$ depending solely on
$d_n,d_s$ and all $N\ge c_0(d_n,d_s)$, the $\alpha$-resolution $\rho
_*(\alpha)$ in the problem $(N_2)$ satisfies:
\begin{enumerate}[(ii)]
\item[(i)] in the range $0\le\epsilon_s\le
c_1(d_n,d_s)N^{-d_s-1/2}\sqrt{\ln(1/\alpha)}$,
\begin{eqnarray*}
\rho_*(\alpha)\geq c_2(d_n,d_s)\sqrt{
\ln(1/\alpha)/N};
\end{eqnarray*}
\item[(ii)] in the range
%
\begin{eqnarray}
\label{range1} &&c_3(d_n,d_s)N^{-d_s-1/2}
\sqrt{\ln(1/\alpha)}\leq\epsilon_s\leq c_1(d_n,d_s)
\sqrt{\ln(1/\alpha)},
\\
&&\rho_*(\alpha) \geq c_4(d_n,d_s)
\bigl[{\epsilon _s} N^{{d_s}+1/2} \bigl[\ln(1/\alpha)
\bigr]^{-1/2} \bigr]^{1/(2d_s+1)}\sqrt{N^{-1}\ln(1/\alpha)}
\nonumber
\\[-8pt]
\\[-8pt]
&&\hphantom{\rho_*(\alpha)}\geq c_5(d_n,d_s)
\epsilon_s^{1/(2d_s+1)} \bigl(\ln(1/\alpha ) \bigr)^{d_s/(2d_s+1)}.
\nonumber
\end{eqnarray}
\end{enumerate}
\end{proposition}

The results of items (i) in Propositions \ref{lowbndC} and \ref
{lowbndD} say that when $d_n,d_s$ are fixed, $N$ is large, and
${\epsilon_n}, {\epsilon_s}$ are small enough so that the problem
parameters are in the parametric zone (i.e., (\ref{epsilonsaresmall}) holds),
Theorem~\ref{themain11} describes ``nearly correctly'' (i.e., up to
factors depending solely on $d_n,d_s,N$ and logarithmic in $N$) the
$\alpha$-resolution in problems $(N_1)$ and $(N_2)$: within such a
factor, the $\alpha$-resolution for problems $(N_1)$, $(N_2)$, same as
for problems $(P_1)$, $(P_2)$, is $\sqrt{\ln(1/\alpha)/N}$. Besides,
items (ii) in Propositions \ref{lowbndC} and \ref{lowbndD} show that
when $(\epsilon_n,\epsilon_s)$ goes ``far beyond'' the range (\ref
{epsilonsaresmall}), the $\alpha$-resolution in problems $(N_1)$,
$(N_2)$ becomes ``much worse'' than $\sqrt{\ln(1/\alpha)/N}$.

\section{Numerical results}\label{sectNumres}
Below we report on {some} numerical experiments with the basic test.
\subsection{Problem $(N_1)$}
The goal of the first series of {simulations} was to quantify
``practical performance'' of the basic test as applied to problem $(N_1)$.
\subsubsection*{Organization of experiments} We consider problem $(N_1)$
on the time horizon $0\leq t<N $ for $N\in\{128,512,1024\}$ with
reliability threshold $\alpha=0.01$. In these {simulations} $d_n=4$,
the frequencies in $\overline{\mathbf{w}}$ are selected at random,
$\epsilon_n=0.01$, and $d_s=4$ (note that $N$ and $\epsilon_n$ are
deliberately chosen \textsl{not} to satisfy~(\ref{epsilonsaresmall})).
As explained in Section~\ref{sectTest}, the above setup specifies the
basic test for problem $(N_1)$, and our goal is to find the ``empirical
resolution'' of this test. To this end, we ran 10 experiments as
follows. In a particular experiment:
\begin{itemize}
\item We draw at random ${\mathbf{w}}\in\cS_4$, a \textsl{shift}
$\bar{s}\in\cS[{\mathbf{w}}]$ and \textsl{basic nuisance} $u\in
\cS^{\epsilon_n}[\overline{\mathbf{w}}]$.
\item We generate a ``true signal'' $x$ according to $x_\lambda
=\lambda\bar{s}+u$, where $\lambda>0$ is (nearly) as small as
possible under the restriction that with $x=x_\lambda$, the basic test
``rejects reliably'' the hypothesis $H_0$, namely, rejects it in every
one of 15 trials with $x=x_\lambda$ and different realizations of the
observation noise $\xi$, see (\ref{eqobs}).\footnote{Since a run of
the test requires solving a nontrivial convex program, it would be too
time-consuming to replace 15 trials with few hundreds of them required
to check reliably that the probability to reject $H_0$, the signal
being $x_\lambda$, is at least the desired $1-\alpha=0.99$.}
\item For the resulting $\lambda$, we compute $\rho=\min_{u\in\cS
^{\epsilon_n}[\overline{\mathbf{w}}]}\|x_\lambda-u\|_\infty$,
which is the output of the experiment. We believe that the collection
of 10 outputs of this type gives a good impression on the ``true
resolution'' of the basic test. As a byproduct of an experiment, we get
also the $\|\cdot\|_\infty$-closest to $x_\lambda$ point $u_x\in\cS
^{\epsilon_n}[\overline{\mathbf{w}}]$; the quantity $r=\|x_\lambda
-u_x\|_2/\sqrt{N}$ can be thought of as a natural in our context
``signal-to-noise ratio.''
\end{itemize}

\textit{The results} are presented in Table~\ref{table1}.
%
\begin{table}
\tablewidth=\textwidth
\tabcolsep=0pt
\caption{Problem $(N_1)$ with $d=d_s=4$, $\alpha
=\epsilon=0.01$}\label{table1} %
\begin{tabular*}{\textwidth}{@{\extracolsep{\fill}}lllllllllllll@{}}
\hline
&\multicolumn{10}{l}{Experiment \#}&&\\[-5pt]
&\multicolumn{10}{l}{\hrulefill}&&\\
&1&2&3&4&5&6&7&8&9&10&Mean&Mean$\times N^{1/2}$\\
\hline
\multicolumn{13}{l}{Experiments with $N=128$:}\\
\quad Resolution&1.10&1.58&1.52&1.51&2.28&1.85&1.12&1.92&1.12&1.82&1.58&17.9\\
\quad Signal/noise&0.70&1.09&0.94&0.95&1.36&1.03&0.77&1.06&0.69&1.10&0.97&11.0\\
[3pt]
\multicolumn{13}{l}{Experiments with $N=512$:}\\
\quad Resolution&0.79&1.30&1.31&0.79&0.79&0.79&0.48&0.79&0.81&0.48&0.83&18.8\\
\quad Signal/noise&0.44&0.71&0.74&0.43&0.43&0.44&0.30&0.50&0.44&0.34&0.48&10.8\\
[3pt]
\multicolumn{13}{l}{Experiments with $N=1024$:}\\
\quad Resolution&0.60&0.92&0.36&0.58&0.46&0.35&0.36&0.56&0.42&0.59&0.52&16.6\\
\quad Signal/noise&0.32&0.56&0.24&0.41&0.31&0.27&0.26&0.43&0.26&0.39&0.35&11.0\\
\hline
\end{tabular*}
\end{table}

%
\begin{figure}
%

\includegraphics{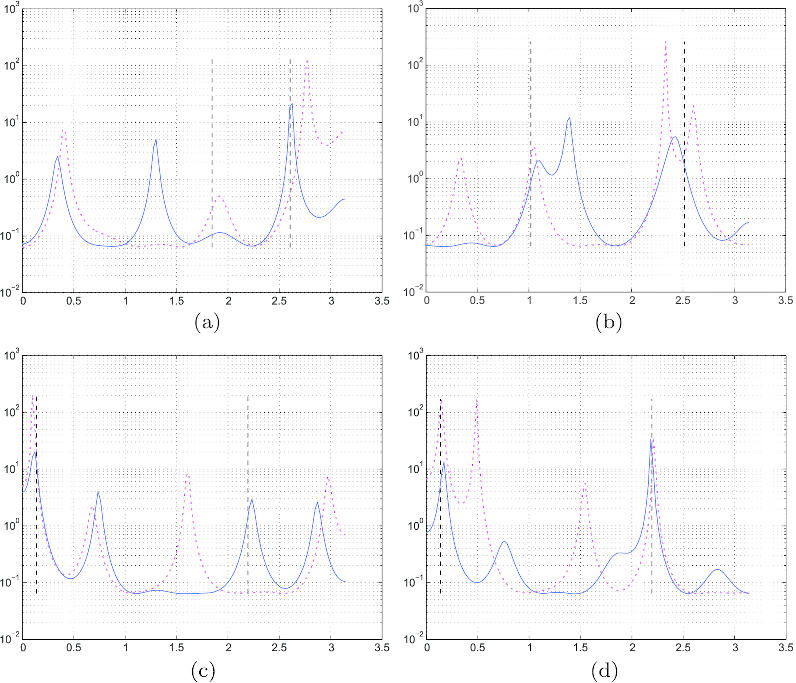}

\caption{MUSIC pseudospectra as built by MATLAB function
\texttt{pmusic}($\cdot$,8).$^\dag$ Dot (magenta): signal plus
nuisance; solid (blue): pure nuisance; dash vertical bars: nuisance
frequencies ($d=4$ elements in $\overline{\mathbf{w}}$ correspond to
2 ``actual'' frequencies). (a), (b) $N=128$; (c), (d) $N=1024$.
$^\dag$In the present setup, $\overline{\mathbf{w}}\cup{\mathbf
{w}}=\{\pm w_j,1\leq j\leq4\}$, which requires the \texttt{pmusic}
parameter \texttt{p} to be set to 8.}\vspace*{-6pt}\label{fig1}
\end{figure}

\noindent We would qualify them as quite compatible with the theory we have
developed: both empirical resolution and empirical signal-to-noise
ratio decreases with $N$ as $N^{-1/2}$. The ``empirically observed''
resolution $\rho$ for which the basic test $(1-\alpha)$-reliably,
$\alpha=0.01$, distinguishes between the hypotheses $H_0$ and
$H_1(\rho)$ associated with problem $(N_1)$ is $\approx6\sqrt{\ln
(N/\alpha)/N}$.

\subsubsection*{Comparison with MUSIC} An evident alternative to the basic
test is (a) to apply the standard MUSIC algorithm \cite{mus1} in order
to recover the spectrum of the observed signal, (b) to delete from this
spectrum the ``nuisance frequencies,'' and (c) to decide from the
remaining data if the signal of interest is present. Our related
numerical results are, to the best of our understanding, strongly in
favor of the basic test. Let us look at Figure~\ref{fig1} where we
present four MUSIC pseudospectra (we use \texttt{pmusic} function from
MATLAB Signal Processing Toolbox) of the observations associated with
signals $x$ obeying the hypothesis $H_1(\rho)$ (magenta) and of the
observations coming from the $\|\cdot\|_\infty$-closest to $x$
nuisance (i.e., obeying the hypotheses $H_0$) $u_x$ (blue). $\rho$ was
chosen large enough for the basic test to accept reliably the
hypothesis $H_1(\rho)$ when it is true. We see that while sometimes
MUSIC pseudospectrum indeed allows to understand which one of the
hypotheses takes place (as it is the case in the example (d)), ``MUSIC
abilities'' in our context are rather limited.\footnote{It should be
noted that MUSIC is designed for a problem different from (and more
complex than) the detection we are interested in, and thus its weakness
relative to a dedicated detection test does not harm algorithm's
well-established reputation.} For example, it is hard to imagine a
routine which would attribute magenta curves in the examples (a)--(c) to
signals, and the blue curves -- to the nuisances.

%
\subsection{Comparison with energy test}
Our objective here is to compare the resolution of the basic test to
that of the test which implements the straightforward idea of how to
discover if the signal $x$ underlying observations (\ref{eqobs}) does
not belong to a {known} nuisance set $\cU\subset\cS$. The test in
question, which we refer to as \textsl{energy test}, is as follows:
given a tolerance $\alpha$ and an observation $y$, we solve the
optimization problem
\[
\Opt(y)=\inf_{u\in\cU} \bigl\|y-u_0^{N-1}
\bigr\|_2^2
\]
and compare the optimal value with the $(1-\alpha)$-quantile
\[
p_N(\alpha)\dvtx \Prob_{\xi\sim\cN(0,I_N)} \bigl\{\|\xi
\|_2^2>p_N( \alpha) \bigr\} =\alpha
\]
of the $\chi^2$-distribution with $N$ degrees of freedom. If $\Opt
(y)>p_N(\alpha)$, we reject the nuisance hypothesis $H_0$ stating that
$x\in\cU$, otherwise we accept the hypothesis. Note that\vspace*{1pt} the basic
test is of a completely similar structure, with $\|F_N(y-u_0^{N-1})\|
_\infty^2$ in the role of $\|y-u_0^{N-1}\|_2^2$ and $q_N^2(\alpha)$
in the role of $p_N(\alpha)$. It is clear that the energy test rejects
$H_0$ when the hypothesis is true with probability at most $\alpha$
(cf.\ item 1$^0$ in Section~\ref{proofThemain1}). In order to {simplify
the presentation}, we restrict this test comparison to the simplest
case of $\cU=\{0\}$, i.e., the case of problem $(P_1)$. Let us start
with some theoretical analysis. Given a natural $d_s>0$ and a real
$\rho>0$, consider the signal hypothesis $\overline{H}_1(\rho)$
stating that
the signal $x$ underlying observations (\ref{eqobs}) satisfies $\|
x_0^{N-1}\|_\infty\geq\rho$ and that $x_t$ is a real algebraic
polynomial of degree $\leq d-1$ of $t\in\bZ$, meaning that $x\in\cS
[\overbrace{0,\ldots,0}^d]$. Observe that with our $\cU=\{0\}$, $\Opt
(y)$ is nothing but
\[
\|y\|_2^2=\bigl\|x_0^{N-1}+\xi
\bigr\|_2^2=\|\xi\|_2^2+2
\xi^Tx_0^{N-1}+\bigl\| x_0^{N-1}
\bigr\|_2^2.
\]
It follows that the hypothesis $H_0$ is accepted whenever the event
\[
\|\xi\|_2^2+2\xi^Tx_0^{N-1}+
\bigl\|x_0^{N-1}\bigr\|_2^2\leq
p_N(\alpha)
\]
takes place. Now, from the standard results on the $\chi^2$
distribution it follows that for every $\alpha\in(0,1)$, for all
large enough values of $N$ with properly chosen absolute constants it holds
\[
\Prob_{\xi\sim\cN(0,I_N)} \bigl\{p_N(\alpha)-\|\xi\|_2^2
\geq \mathrm{O}(1)\sqrt {N\ln(1/\alpha)} \bigr\}\geq \mathrm{O}(1),
\]
whence also
\[
\Prob_{\xi\sim\cN(0,I_N)} \bigl\{ \bigl\{p_N(\alpha)-\|\xi
\|_2^2\geq \mathrm{O}(1)\sqrt{N\ln(1/\alpha)} \bigr\}\cap
\bigl\{ \xi^Tx_0^{N-1}\leq0 \bigr\} \bigr\} \geq
\mathrm{O}(1).
\]
As a result, whenever $x\in\cS$ satisfies $\|x_0^{N-1}\|_2^2\leq
\mathrm{O}(1)\sqrt{N\ln(1/\alpha)}$, the probability to accept $H_0$, the
true signal being $x$, is at least $\mathrm{O}(1)$, provided that $N$ is large
enough. On the other hand, for a given $d_s$ and large $N$ there exists
a polynomial $x$ of degree $d_s-1$ such that $\|x_0^{N-1}\|_2\leq
d_s^{-1}N^{1/2}\|x_0^{N-1}\|_\infty$, see the proof of Proposition~\ref{lowbndAB}. It immediately follows that with $d_s\geq1$ and
(small enough) $\alpha>0$ fixed, the energy test cannot distinguish
$(1-\alpha)$-reliably between the hypotheses $H_0$ and $\overline
{H}_1(\rho)$, provided that
%
\begin{equation}
\label{lowbndEnergy} \rho=\mathrm{O}(1)d_s \bigl[N^{-1}\ln(1/
\alpha) \bigr]^{1/4}
\end{equation}
and $N$ is large enough. In other words, with $d_s$ and (small enough)
$\alpha$ fixed, the resolution of the energy test in problem $(P_1)$
admits, for large $N$, the lower bound (\ref{lowbndEnergy}). Note that
as $N$ grows, this bound goes to 0 as $N^{-1/4}$, while the resolution
of the basic test goes to 0 as $N^{-1/2}\sqrt{\ln(N)}$ (Theorem~\ref
{themain1}).
We conclude that the basic test provably outperforms the energy test as
$N\to\infty$.
The goal of the {simulation} experiments we are about to report is to
investigate this phenomenon numerically.

\subsubsection*{Organization of experiments} In the {simulations} to
follow, the basic test and the energy test were tuned to
0.99-reliability ($\alpha=0.01$) and used on time horizons $N\in\{
256,1024,4096\}$. For a fixed $N$, and every value of the ``resolution
parameter'' $\rho$ from the equidistant grid on $[0,4]$ with {the grid
step of} $0.05$, we run 10\,000 {simulations} as follows:
\begin{itemize}
\item we generate $z\in\cS[{\mathbf{w}}]$, and specify signal $x$ as
$\rho z/\|z_0^{N-1}\|_\infty$;
\item we generate $y$ according to (\ref{eqobs}) and run on the
observations $y$ the basic test and the energy test.
\end{itemize}
For {each test,} the outcome of a series of 10\,000 {simulations} is the
empirical probability $p$ of rejecting the nuisance hypothesis $H_0$
(which states that the signal underlying the observations is
identically zero). For $\rho=0$, $p$ is the (empirical) probability of
false alarm (rejecting $H_0$ when it is true), and we want it to be
small (about $\alpha=0.01$). For $\rho>0$, $p$ is the empirical
probability of successful detection of an actually present signal, and
we want it to be close to 1 (about $1-\alpha=0.99$). Given that $p\leq
\alpha$ when $\rho=0$, the performance of a test can be quantified as
the smallest value $\rho_*$ of $\rho$ for which $p$ is at least
$1-\alpha$ (the less is $\rho_*$, the better).

We use 4-element collections $\mathbf{w}$ (i.e., $d_s=4$), and for
every $N$ and $\rho$ run two 10\,000-element series of {simulations}
differing in how we select $\mathbf{w}$ and $z$. In the first series
(``random signals''), $\mathbf{w}$ is selected at random, and $z$ is a
random combination of the corresponding harmonic oscillations. In the
second series (``bad signal'') we use ${\mathbf{w}}=\{0,0,0,0\}$, and
$z$ is the algebraic polynomial of degree 3 with the largest, among
these polynomials, ratio of $\|z_0^{N_1}\|_\infty/\|z_0^{N-1}\|_2$. In
the latter case, only the realisation of noise varied from one
experiment to another.

\textit{The results} of our experiments are presented in Table~\ref{table2}. They are in full accordance to what is suggested by our
theoretical analysis; for $N=256$, both tests exhibit nearly the same
empirical performance. As $N$ grows, the empirical performances of both
tests improve, and the ``performance gap'' (which, as expected, is in
favor of the basic test) grows.
%
\begin{sidewaystable}
\tablewidth=\textwidth
\tabcolsep=0pt
\caption{Basic test vs. Energy test, problem $(P_1)$
with $d_s=4$. $p(\mathrm{B})$, $p(\mathrm{E})$: empirical probabilities,
taken over 10\,000 trials, of detecting signal using the basic test (B)
and the energy test (E). $\rho_*(\cdot)$ is the smallest $\rho$ for
which $p(\cdot)\geq1-\alpha=0.99$}\label{table2}
\vspace*{6pt}
\begin{tabular*}{\textwidth}{@{\extracolsep{\fill}}lllllllllllllll@{}}
\multicolumn{15}{l}{$N=256$, random signals ($\rho_*(\mathrm
{B})\approx1.10,\rho_*(\mathrm{E})\approx1.35$):}\\
\hline
& \multicolumn{14}{l}{$\rho$}\\[-5pt]
&\multicolumn{14}{l}{\hrulefill}\\
&0.00&0.95&1.00&1.05&1.10&1.15&1.20&1.25&1.30&1.35&1.40&1.45&1.50&1.55\\
\hline
$p(\mathrm
{B})$&0.010&0.960&0.977&0.987&0.993&0.998&0.999&0.999&1.000&1.000&1.000 &1.000&1.000&1.000 \\
$p(\mathrm
{E})$&0.011&0.710&0.779&0.842&0.887&0.933&0.956&0.974&0.987&0.994&0.997&0.999&1.000&1.000\\
\hline
\end{tabular*}
\begin{tabular*}{\textwidth}{@{\extracolsep{\fill}}lllllllllllllll@{}}
\\[-5pt]
\multicolumn{15}{l}{$N=256$, ``bad'' signal ($\rho_*(\mathrm
{B})\approx2.65,\rho_*(\mathrm{E})\approx2.75$):}\\
\hline
& \multicolumn{14}{l}{$\rho$}\\[-5pt]
&\multicolumn{14}{l}{\hrulefill}\\
&0.00&2.50&2.55&2.60&2.65&2.70&2.75&2.80&2.85&2.90&2.95&3.00&3.05&3.10\\
\hline
$p(\mathrm
{B})$&0.010&0.973&0.984&0.987&0.991&0.994&0.998&0.998&0.999&1.000&1.000&1.000&1.000&1.000\\
$p(\mathrm
{E})$&0.011&0.941&0.956&0.971&0.978&0.984&0.990&0.993&0.995&0.997&0.998&0.998&1.000&1.000\\
\hline
\end{tabular*}
\begin{tabular*}{\textwidth}{@{\extracolsep{\fill}}lllllllllllllll@{}}
\\[-5pt]
\multicolumn{15}{l}{$N=1024$, random signals ($\rho_*(\mathrm
{B})\approx0.60,\rho_*(\mathrm{E})\approx0.90$):}\\
\hline
& \multicolumn{14}{l}{$\rho$}\\[-5pt]
&\multicolumn{14}{l}{\hrulefill}\\
&0.00&0.50&0.55&0.60&0.65&0.70&0.75&0.80&0.85&0.90&0.95&1.00&1.05&1.10\\
\hline
$p(\mathrm
{B})$&0.010&0.960&0.986&0.997&1.000&1.000&1.000&1.000&1.000&1.000&1.000&1.000&1.000&1.000\\
$p(\mathrm
{E})$&0.010&0.303&0.421&0.559&0.686&0.795&0.886&0.938&0.974&0.992&0.997&0.999&1.000&1.000\\
\hline
\end{tabular*}
\end{sidewaystable}
\setcounter{table}{1}
\begin{sidewaystable}
\tablewidth=\textwidth
\tabcolsep=0pt
\caption{\textit{Continued}}
\vspace*{6pt}
\begin{tabular*}{\textwidth}{@{\extracolsep{\fill}}lllllllllllllllll@{}}
\\[-5pt]
\multicolumn{17}{l}{$N=1024$, ``bad'' signal ($\rho_*(\mathrm
{B})\approx1.40,\rho_*(\mathrm{E})\approx1.90$):}
\\\hline
& \multicolumn{16}{l}{$\rho$}\\[-5pt]
&\multicolumn{16}{l}{\hrulefill}\\
&0.00&1.30&1.35&1.40&1.45&1.50&1.55&1.60&1.65&1.70&1.75&1.80&1.85&1.90&1.95&2.00\\
\hline
$p(\mathrm
{B})$&0.011&0.960&0.980&0.990&0.997&0.998&0.999&1.000&1.000&1.000&1.000&1.000&1.000&1.000&1.000&1.000\\
$p(\mathrm
{E})$&0.011&0.505&0.564&0.633&0.703&0.770&0.819&0.863&0.903&0.932&0.960&0.971&0.987&0.993&0.996&0.997\\
\hline
\end{tabular*}
\begin{tabular*}{\textwidth}{@{\extracolsep{\fill}}llllllllllll@{}}
\\[-5pt]
\multicolumn{12}{l}{$N=4096$, random signals ($\rho_*(\mathrm
{B})\approx0.30,\rho_*(\mathrm{E})\approx0.65$):}\\
\hline
& \multicolumn{11}{l}{$\rho$}\\[-5pt]
&\multicolumn{11}{l}{\hrulefill}\\
&0.00&0.25&0.30&0.35&0.40&0.45&0.50&0.55&0.60&0.65&0.70\\
\hline
$p(\mathrm
{B})$&0.009&0.931&0.993&0.999&1.000&1.000&1.000&1.000&1.000&1.000&1.000\\
$p(\mathrm
{E})$&0.009&0.084&0.165&0.291&0.472&0.667&0.823&0.930&0.980&0.997&1.000\\
\hline
\end{tabular*}
\begin{tabular*}{\textwidth}{@{\extracolsep{\fill}}llllllllllllllllll@{}}
\\[-5pt]
\multicolumn{18}{l}{$N=4096$, ``bad'' signal ($\rho_*(\mathrm
{B})\approx0.75,\rho_*(\mathrm{E})\approx1.35$):}\\
\hline
& \multicolumn{17}{l}{$\rho$}\\[-5pt]
&\multicolumn{17}{l}{\hrulefill}\\
&0.00&0.70&0.75&0.80&0.85&0.90&0.95&1.00&1.05&1.10&1.15&1.20&1.25&1.30&1.35&1.40&1.45\\
\hline
$p(\mathrm
{B})$&0.010&0.975&0.994&0.999&1.000&1.000&1.000&1.000&1.000&1.000&1.000&1.000&1.000&1.000&1.000&1.000&1.000\\
$p(\mathrm
{E})$&0.012&0.184&0.223&0.290&0.376&0.477&0.567&0.676&0.767&0.843&0.899&0.945&0.975&0.986&0.995&0.999&1.000\\
\hline
\end{tabular*}
\end{sidewaystable}

\section{Proofs}\label{sec:proofs}
\subsection{Preliminaries}\label{sec:outline}
\subsubsection*{Notation} In what follows, for a {complex valued}
polynomial $p(\zeta)=\sum_{k=0}^mp_k\zeta^k$, we denote
\[
\bigl\|p(\cdot)\bigr\|_\infty=\max_{\zeta\in\C,|\zeta|=1} \bigl|p(\zeta)\bigr|
\]
and denote by
\[
|p|_{s}=\bigl\|[p_0;p_1;\ldots ;p_m]
\bigr\|_s,\qquad 1\leq s\leq\infty,
\]
the $\ell_s$-norm of the vector of coefficients, so that
\[
\bigl\|p(\cdot)\bigr\|_2^2:={1\over2\uppi}
\oint_{|\zeta|=1} \bigl|p(\zeta )\bigr|^2|\mathrm{d}\zeta|=|p|_2^2.
\]

\textit{The key fact} underlying Theorems \ref{themain1}, \ref
{themain11} is the following
proposition.

\begin{proposition}\label{MainProp} Let $d$, $N$ be positive integers
and $s\in\cS_d$. Then
%
\begin{equation}
\label{MainPropRel} \|F_Ns\|{_\infty}\geq c(d) N^{1/2}
\bigl\|s_0^{N-1}\bigr\|_\infty,
\end{equation}
where $c(d)>0$ is a universal nonincreasing function of $d$. One can take
%
\begin{equation}
\label{cofd} c(d)=\mathrm{O}(1)/ \bigl(d^3\ln(2d) \bigr)
\end{equation}
with properly selected positive absolute constant $\mathrm{O}(1)$.
\end{proposition}

\begin{pf} Let us fix $s\in\cS_d$; we intend to prove that $s$
obeys (\ref{MainPropRel}). Let ${\mathbf{u}}=\{\omega_1,\ldots,\omega
_d\}$ be a symmetric $\mod  2\uppi$ collection such that $s\in\cS
[{\mathbf{u}}]$, and let
\[
p_{\mathbf{u}}(\zeta) = \prod_{\ell=1}^d
\bigl(1-\exp\{\imath\omega _\ell\}\zeta \bigr),
\]
so that $p_{\mathbf{u}}(\Delta)s\equiv0$. {Further, let} $M$ be the
index of the largest in magnitude of the reals $s_0,s_1,\ldots,s_{N-1}$,
so that
%
\begin{equation}
|s_M|=\bigl\|s_0^{N-1}\bigr\|_\infty.
\end{equation}
We can w.l.o.g. assume that $M \geq(N-1)/2$. Indeed, otherwise we
could pass from $s$ to the ``reversed'' sequence $s'\in\cS_d$:
$s^\prime_t=s_{N-t-1}$, $t\in\Z$, which would not affect the
validity of our target relation (\ref{MainPropRel}) and would convert
$M<(N-1)/2$ into $M'=N-1-M\geq(N-1)/2$.

1$^0$. We need the following technical
lemma.

\begin{lemma}\label{MainLemma} Let $d$ be a positive integer, and let
${\mathbf{u}}=\{\upsilon_1,\ldots,\upsilon_d\}\in\Omega_d$. For every
integer $m$ satisfying
%
\begin{equation}
\label{mislargeenough} m\geq m(d):=d\Ceil \biggl(5d\max \biggl[2, \frac{1}{2}\ln
(2d) \biggr] \biggr)
\end{equation}
one can point out real polynomials $q(\zeta)=\sum_{j=1}^mq_j\zeta^j$
and $r(\zeta)=1+\sum_{j=1}^{m-d}r_j\zeta^j$ such that
%
\begin{eqnarray}
\label{divisible} 1-q(\zeta)=p_{\mathbf{u}}(\zeta)r(\zeta),
\end{eqnarray}
and
%
\begin{equation}
\label{propofq} |q|_2\leq C_1(d)/\sqrt{m},\qquad
\mbox{where } C_1(d)= 3\e d^{3/2}\sqrt{\ln(2d)}.
\end{equation}
\end{lemma}

The proof of Lemma~\ref{MainLemma} is presented in the \hyperref[app]{Appendix}.

2$^0$. The following statement is immediate:

\begin{lemma}\label{L1F}
Let $m\le{(N-1)/2}$, and $g\in\C^N$ be such that $g_i=0$ for $i>m$.
Let $h\in\C^N$ be the discrete autoconvolution of $g$, that is, the
vector with entries
$h_k=\sum_{0\le i,j\le m, i+j=k} g_ig_j$, $0\leq k\leq2m$ and with
zero remaining $N-2m-1$ entries. Then
\[
\|F_N h\|_1= \sqrt{N} \|g\|_2^2.
\]
\end{lemma}

\begin{pf} We have
\begin{eqnarray*}
[F_Nh](\mu)&=& N^{-1/2}\sum_{0\leq t\leq2m}
\biggl[\sum_{0\leq j,k\leq m,
j+k=t}g_jg_k
\biggr]\mu^t
\\
&=&N^{-1/2}\sum_{0\leq t\leq2m}\sum
_{0\leq j,k\leq m,j+k=t} \bigl(g_j\mu ^{j} \bigr)
\bigl(g_k\mu^{k} \bigr)
\\
&=&N^{1/2} \Biggl[N^{-1/2}\sum_{j=0}^{m}g_j
\mu^{j} \Biggr]^2.
\end{eqnarray*} %
Invoking the Parseval identity, we conclude that
\[
\|F_Nh\|_1=N^{1/2}\|F_Ng
\|^2_2=N^{1/2}\|g\|^2_2.
\]
\upqed
\end{pf}

3$^0$. Let
%
\begin{equation}
\label{let1} N>60d^2\ln(2d),
\end{equation}
and let
%
\begin{equation}
\label{miz} m=\Floor \biggl({N-1\over4} \biggr).
\end{equation}
Then (\ref{mislargeenough}) is satisfied, and, according to Lemma~\ref
{MainLemma}, there exists a polynomial $q(\zeta)=\sum_{j=1}^mq_j\zeta
^j$ such that
\[
1-q(\zeta)=p_{\mathbf{u}}(\zeta)r(\zeta),\qquad |q|_2\leq
C_1(d)/\sqrt{m},
\]
with some polynomial $r$. Setting $q^+(\zeta)=q^2(\zeta)=\sum_{j=1}^{2m}q^+_j\zeta^j$, we get
\begin{eqnarray*}
q^+_k&=&\sum_{1\le i,j\le m, i+j=k}q_iq_j,
\\
\bigl(1-q^+(\Delta) \bigr)s&=& \bigl(1+q(\Delta) \bigr) \bigl(1-q(\Delta) \bigr)s=
\bigl(1+q(\Delta ) \bigr)r(\Delta)p_{\mathbf{u}}(\Delta)s\equiv0,
\end{eqnarray*}
whence
$s_M=\sum_{i=1}^{2m}s_{M-i}q^+_i$ (note that $M\geq(N-1)/2\geq2m$ by
(\ref{miz})). Let now $h\in\R^N$ be the vector with coordinates
\[
h_i=\lleft\{ %
\begin{array} {l@{\qquad}l}
q^+_{M-i}, &i=M-1,\ldots,M-2m,
\\
0,& \mbox{otherwise}. \end{array} %
\rright.
\]
Note that by Lemma~\ref{L1F} and due to $|q|_2\leq C_1(d)/\sqrt{m}$
one has
\[
\|F_Nh\|_1\le N^{1/2}|q|_2^2
\le C_1^2(d)\sqrt{N}/m.
\]
We have
\[
\bigl\|s_{0}^{N-1}\bigr\|_{\infty}=|s_M|= \Biggl
\llvert \sum_{i=1}^{2m} q^+_is_{M-i}
\Biggr\rrvert =\bigl| \bigl\langle h,s_{0}^{N-1} \bigr\rangle\bigr|=\bigl|
\langle F_Nh,F_Ns\rangle\bigr|\le\|F_Nh
\|_1\|F_Ns\|_\infty,
\]
where the last equality is given by the fact that $F_N$ is unitary, whence
\[
\|F_Ns\|_\infty\geq\bigl\|s_0^{N-1}
\bigr\|_\infty/\|F_Nh\|_1\geq{m\over
C_1^2(d)\sqrt{N}}
\bigl\|s_0^{N-1}\bigr\|_\infty.
\]
Invoking (\ref{miz}), (\ref{let1}), and (\ref{propofq}), we see that
for $N$ satisfying (\ref{let1}) our target relation (\ref
{MainPropRel}) indeed holds true, provided that
%
\begin{equation}
\label{providedthat} c(d)\leq \mathrm{O}(1) \bigl[d^3\ln(2d)
\bigr]^{-1}
\end{equation}
with properly selected positive absolute constant $\mathrm{O}(1)$.

4$^0$. It remains to verify (\ref{MainPropRel}) when
$N\leq60 d^2\ln(2d)$. Since $F_N$ is unitary, we have $\|s_0^{N-1}\|
_\infty\leq\|s_0^{N-1}\|_2=\|F_Ns\|_2\leq\|F_N s\|_\infty\sqrt {N}$, whence
\[
\|F_N s\|_\infty\geq N^{-1/2}\bigl\|s_0^{N-1}
\bigr\|_\infty\geq N^{-1} \bigl[N^{1/2}\bigl\|
s_0^{N-1}\bigr\|_\infty \bigr]\geq{1\over60 d^2\ln(2d)}
\bigl[N^{1/2}\bigl\|s_0^{N-1}\bigr\| _\infty \bigr],
\]
which completes the proof.
\end{pf}

\subsection{Proof of Theorem \texorpdfstring{\protect\ref{themain1}}{3.1}}\label
{sec:themain1}\label{proofThemain1}

1$^0$. Let us prove the result for the basic test,
let it be denoted $\widehat{\varphi}$, as applied to the problem
$(P_2)$; note that $(P_1)$ is the particular case of $(P_2)$
corresponding to $d_n=0$.

We have to show that under the premise of the theorem $\varepsilon
_0(\widehat{\varphi}; H_0)\le\alpha$ and $ \varepsilon_1(\widehat
{\varphi}; H_{1}(\rho))\le\alpha$. The first bound is evident.
Indeed, let $\Xi_\alpha=\{\xi\dvtx  \|F_N \xi\|_\infty\leq q_N(\alpha
)\}$, so\vspace*{1pt} that $\Prob_{\xi\sim\cN(0,I_N)}\{\xi\in\Xi_\alpha\}
\geq1-\alpha$. Under the hypothesis $H_0$, the set $\mathcal{Z}$
from (\ref{convexproblem1}) contains the true signal $x_0^{N-1}$, so
that the optimal value $\Opt_\mathcal{Z}(y)$ in (\ref{convexproblem0})
is at most $\|F_N \xi\|_\infty$. It follows that when $\xi\in\Xi
_\alpha$ (which happens with probability $\geq1-\alpha$) we have
${\Opt_\mathcal{Z}(y)}\leq q_N(\alpha)$, and the basic test will
therefore accept $H_0$. We conclude that $\varepsilon_0(\widehat
{\varphi}; H_0)\le\Prob\{\xi\notin\Xi_\alpha\}\leq\alpha$.

2$^0$. Now let $x\in H_1(\rho)$, that is,
$x=s+u$, where $s\in\cS[\mathbf{w}]$ for some $\mathbf{w}\in\Omega
_{d_s}$, $u\in\cS[\overline{\mathbf{w}}]$, and
\[
\bigl\|[x-z]_0^{N-1}\bigr\|_\infty\geq\rho\qquad \forall z\in
\cS[{ \overline {\mathbf{w}}}].
\]
Let $z\in\cS[{\overline{\mathbf{w}}}]$, and let $s=x-z$. Then $s\in
\cS_{d_*}$, $d_*=d_n+d_s$, and $\|s_0^{N-1}\|_\infty\geq\rho$,
whence, by Proposition~\ref{MainProp},
\[
\|F_Ns\|_\infty\geq c(d_*)N^{1/2}
\bigl\|s_0^{N-1}\bigr\|_\infty\geq c(d_*)N^{1/2}\rho.
\]
It follows that the optimal value $\Opt_\mathcal{Z}(y)$ in (\ref
{convexproblem0}) is at least $c(d_*)N^{1/2}\rho-\|F_N\xi\|_\infty$.
Recalling the definition of $q_N(\alpha)$, we conclude that
\[
\Prob \bigl\{\Opt_\mathcal{Z}(y)>c(d_*)N^{1/2}
\rho-q_N(\alpha) \bigr\}\geq 1-\alpha
\]
as soon as
%
\begin{equation}
\label{when} \rho>{{2q_N(\alpha)\over c(d_*)\sqrt{N}}},
\end{equation}
and the probability to reject $H_1(\rho)$ when the hypothesis is true
is $\leq\alpha$. {We} see that {for the proper choice of the absolute
constant $\mathrm{O}(1)$} under the premise of Theorem~\ref{themain1} one has
\mbox{$\varepsilon_1(\widehat{\varphi}, H_1(\rho))\leq\alpha$}.

\subsection{Proof of Theorem \texorpdfstring{\protect\ref{themain11}}{3.2}}

1$^0$. We start with the following simple

\begin{lemma}\label{lemdev} Let $d$ and $N$ be positive integers, let
$\epsilon\geq0$, let ${\mathbf{u}}=\{\upsilon_1,\ldots,\upsilon_d\}$
be a symmetric $\mod 2\uppi$ $d$-element collection of reals.
Whenever ${w}\in\cS^{N,\epsilon}[{\mathbf{u}}]$, there exists a
decomposition $w=s+z$ such that $s\in\cS[{\mathbf{u}}]$ and
%
\begin{equation}
\label{trivbound} \bigl\|z_0^{N-1}\bigr\|_\infty\leq
N^d\epsilon.
\end{equation}
\end{lemma}

\begin{pf} Let $p_{\mathbf{u}}(\zeta)=\prod_{\ell=1}^d(1-\exp
\{\imath\upsilon_\ell\}\zeta)$ and $r=p_{\mathbf{u}}(\Delta)w$,
so that $\|r_0^{N-1}\|_\infty\leq\epsilon$ due to $w\in\cS
^{N,\epsilon}[{\mathbf{u}}]$. Let, further, $\delta$ be the discrete
convolution unit (i.e., $\delta\in\cS$ is given by $\delta_0=1$,
$\delta_t=0$, $t\neq0$). For $\ell=1,\ldots,d$, let $\gamma^{(\ell)}$
be a two-sided complex-valued sequence obtained from the sequence $\{
\exp\{\imath\upsilon_\ell t\}\}_{t\in\Z}$ by replacing the terms
with negative indexes with zeros, and let $r^+$ be obtained by similar
operation from the sequence $r$. Let us set
\[
\chi=\gamma^{(1)}*\gamma^{(2)}*\cdots*\gamma^{(d)}*r^+,
\]
where $*$ stands for discrete time convolution. It is immediately seen
that $\chi$ is a real-valued two-sided sequence which vanishes for
$t<0$ and satisfies the finite-difference equation $p_{\mathbf
{u}}(\Delta)\chi=r^+$ (due to the evident relation $(1-\exp\{\imath
\upsilon_k\}\Delta)\gamma^{(k)}=\delta$). It follows that
$(p_{\mathbf{u}}(\Delta)(w-\chi))_t=0$ for $t=0,1,\ldots\,$, which (along
with the fact that all the roots of $p_{\mathbf{u}}(\zeta)$ are
nonzero) implies that the sequence $s=w-\chi$ can be modified on the
domain $t<0$ so that $p_{\mathbf{u}}(\Delta)s\equiv0$. Then $z=w-s$
coincides with $\chi$ on the domain $t\ge0$, and $w=s+z$ with $s\in
\cS[{\mathbf{u}}]$ and $z_t=\chi_t$, $t=0,1,\ldots\,$. It remains to note
that for two-sided complex-valued sequences $\mu,\nu$ starting at
$t=0$ we clearly have $\|[\mu*\nu]_0^{N-1}\|_\infty\leq\|\mu
_0^{N-1}\|_1\|\nu_0^{N-1}\|_\infty$. Applying this rule recursively
and taking into account that $\|[\gamma^{(\ell)}]_0^{N-1}\|_1=N$, we
get the recurrence
\[
\bigl\| \bigl[\gamma^{(1)}*\cdots*\gamma^{(\ell+1)}*r^+
\bigr]_0^{N-1} \bigr\|_\infty\leq N \bigl\| \bigl[
\gamma^{(1)}*\cdots*\gamma^{(\ell)}*r^+ \bigr]_0^{N-1}
\bigr\|_\infty,\qquad \ell=0,1,\ldots,d-1,
\]
whence $\|\chi_0^{N-1}\|_\infty\leq N^d\epsilon$. Since $\chi
_t=z_t$ for $t=0,1,\ldots\,$, (\ref{trivbound}) follows.
\end{pf}

2$^0$. We are ready to prove Theorem~\ref
{themain11}. It suffices to consider the case of problem $(N_2)$
(problem $(N_1)$ is the particular case of $(N_2)$ corresponding to
$\epsilon_s=0$). The fact that for the basic test $\widehat{\varphi
}$ one has $\varepsilon_0(\widehat{\bar{\phi}})\leq\alpha$ can be
verified exactly as in the case of Theorem~\ref{themain1}. Let us
prove that under the premise of Theorem~\ref{themain11} we have
$\varepsilon_1(\widehat{\varphi})\leq\alpha$ as well. To this end
let the signal $x$ underlying (\ref{eqobs}) belong to $H_1(\rho)$, so
that $x=r+u$ for some $u\in\cS^{N,\epsilon_n}[\overline{\mathbf
{w}}]$ and some $r\in\cS^{N,\epsilon_n}[\mathbf{w}]$ with
$d_n$-element collection $\overline{\mathbf{w}}$ and $d_s$-element
collection $\mathbf{w}$, both symmetric $\mod 2\uppi$. Let also
$z\in\cS^{N,\epsilon_n}[\overline{\mathbf{w}}]$. Since $x\in
H_1(\rho)$, we have
%
\begin{equation}
\label{sincehave} \bigl\|[x-z]_0^{N-1}\bigr\|_\infty\geq\rho.
\end{equation}
Applying Lemma~\ref{lemdev} to $r$, $u$, $z$, we get the decompositions
%
\begin{eqnarray}
\label{getdec} && x=s+s'+v'\dvtx s\in\cS[\mathbf{w}],
s'\in\cS[\overline{\mathbf {w}}], \bigl\| \bigl[v'
\bigr]_0^{N-1}\bigr\|_\infty\leq N^{d_n}
\epsilon_n+N^{d_s}\epsilon _s,
\nonumber
\\
&&z=s''+v''\dvtx
s''\in\cS[\overline{\mathbf{w}}], \bigl\|
\bigl[v'' \bigr]_0^{N-1}\bigr\|
_\infty\leq N^{d_n}\epsilon_n,
\nonumber
\\
&&\quad \Rightarrow\quad w:=x-z=\bar{s}+\bar{v},
\\
&&\hphantom{\quad {}\Rightarrow{}\quad}\bar{s}=s+s'-s''
\in\cS[\overline{ \mathbf{w}}\cup{\mathbf {w}}],
\nonumber
\\
&&\hphantom{\quad {}\Rightarrow{} \quad} \bar{v}=v'-v'',
\bigl\| \bar{v}_0^{N-1}\bigr\|_\infty\leq\sigma
:=2N^{d_n}\epsilon_n+N^{d_s}\epsilon_s.
\nonumber
\end{eqnarray} %
Now, (\ref{sincehave}) implies that $\|w_0^{N-1}\|\geq\rho$, whence,
by (\ref{getdec}),
\[
\bigl\|\bar{s}_0^{N-1}\bigr\|_\infty\geq\widehat{\rho}:=\rho-
\sigma.
\]
Assuming that $\widehat{\rho}>0$, noting that $\bar{s}\in\cS
[\widetilde{\mathbf{w}}]$ for $(d_*=d_n+d_s)$-element symmetric
$\mod 2\uppi$ collection $\widetilde{\mathbf{w}}$ and invoking
Proposition~\ref{MainProp}, we get
\[
\|F_N\bar{s}\|_\infty\geq c(d_*)N^{1/2}\widehat{
\rho}.
\]
Taking into account that $\|F_N\bar{v}\|_\infty\hspace*{-0.1pt}\leq\hspace*{-0.1pt}\|\bar
{v}_0^{N-1}\|_2\hspace*{-0.1pt}\leq\hspace*{-0.1pt} N^{1/2}\|\bar{v}_0^{N-1}\|_\infty$ and (\ref
{getdec}), we get
also $\|F_N\bar{v}\|_\infty\hspace*{-0.1pt}\leq N^{1/2}\sigma$. Combining these
observations, we get
\begin{eqnarray*}
\bigl\|F_N[x-z]\bigr\|_\infty&=& \|F_N
\bar{s}+F_N\bar{v}\|_\infty\geq\| F_N\bar{s}
\|_\infty-\|F_N\bar{v}\|_\infty
\\
& \geq& c(d_*)N^{1/2}\bar {\rho}-N^{1/2}\sigma
=c(d_*)N^{1/2} [\rho-\sigma ]-N^{1/2}\sigma
\\
& =& \underbrace{c(d_*)N^{1/2} \bigl[\rho- \bigl[1+c_*^{-1}(d_*)
\bigr]\sigma \bigr]}_{\vartheta}.
\end{eqnarray*} %
Since $z\in\cS^{N,\epsilon_n}[\overline{\mathbf{w}}]$ is
arbitrary, we conclude that the optimal value $\Opt_\mathcal{Z}(y)$
in (\ref{convexproblem0}) is at least $\vartheta-\|F_N\xi\|_\infty
$, so that
%
\begin{equation}
\label{sothat17} \Prob \bigl\{\Opt_\mathcal{Z}(y)>\vartheta-q_N(
\alpha) \bigr\}\geq1-\alpha.
\end{equation}
It remains to note that with properly selected positive absolute
constants $\mathrm{O}(1)$'s in (\ref{resolution}) and (\ref
{epsilonsaresmall}), these restrictions on $\rho$, $\epsilon_n$,
$\epsilon_s$ ensure that $\vartheta>2q_N(\alpha)$ (see (\ref
{getdec}), (\ref{cofd})), and therefore (\ref{sothat17}) implies the
desired bound $\varepsilon_1(\widehat{\varphi})\leq\alpha$.

\subsection{Proof of Proposition \texorpdfstring{\protect\ref{lowbndAB}}{4.1}}\label{sec:proof40}
Here we prove the lower resolution bound for problem $(P_2)$. The
result of the proposition for the setting $(P_1)$ may be obtained by an
immediate modification of the proof below for $d_n=0$, $p_{\overline
{\mathbf{w}}}(\cdot)= 1$, and $z\equiv0$.
Below we use notation $\kappa_i$ for positive absolute constants.

1$^0$. Note that for every integer ${d_s}>0$ there exists a
real polynomial $q_{d_s}$ on $[0,1]$ of degree ${d_s}-1$ such that
$\int_0^1q^2_{d_s}(r)\,\mathrm{d}r=1$ and $\max_{0\leq r\leq
1}|q_{d_s}(r)|=q_{d_s}(1)={d_s}$.\footnote{Indeed, let a real
polynomial $p$ of degree $\kappa$ satisfy $p(1)=1$ (and thus $\sum_{i=0}^{\kappa} p_i=1$). Then the vector of coefficients $p^*$ of the
polynomial $p^*_{\kappa}$ of the minimal $L_2(0,1)$-norm is the
optimal solution to the convex quadratic optimization problem
\[
\min \Biggl\{\int_0^1 \Biggl(\sum
_{i=0}^\kappa p_it^i
\Biggr)^2\,\mathrm{d}t=\sum_{i,j=0}^\kappa
{p_ip_j\over i+j+1}, \mbox{subject to} \sum_{i=0}^\kappa
p_i=1 \Biggr\}.
\]
The latter problem can be solved explicitly and its optimal value is
$(\kappa+1)^{-2}$.} Let
$q= \{q_{t}={N}^{-1/2}q_{{d_s}}(t/N) \}_{t=-\infty
}^{+\infty}$, $\lambda>0$, and consider the two-sided sequence
\[
\bar{x}= \bigl\{\bar{x}_t=\theta_N\lambda(-1)^tq_t
\bigr\} _{t=-\infty}^{+\infty},
\]
where $\theta_N>0$ is chosen in such a way that
\[
\bigl\|\bar{x}_{0}^{N-1}\bigr\|_2=\lambda.
\]
Note that $\theta_N$ given by this requirement does not depend on
$\lambda$ and that $\theta_N\to1$ as $N\to\infty$ due to $\int_0^1q_{d_s}^2(r)\,\mathrm{d}r=1$.
We have
%
\begin{equation}
\label{eq222} |\bar{x}_{N-1}|/\bigl\|\bar{x}_0^{N-1}
\bigr\|_2\to{|q_{d_s}(1)|\over\sqrt {\int_0^1q_{d_s}^2(r)\,\mathrm{d}r}}={d_s} \qquad \mbox{as } N
\to \infty.
\end{equation}
The derivative $q'_{d_s}(r)$ of the polynomial $q_{d_s}$ satisfies
\[
\max_{0\leq r\leq1}\bigl|q'_{d_s}(r)\bigr|
\leq2(d_s-1)^2\max_{0\leq r\leq
1}\bigl|q_{d_s}(r)\bigr|
\leq2d_s^3
\]
(the first inequality in this chain follows from {Markov brothers'}
inequality). We conclude that for properly selected $\kappa_1\geq1$
and all $N\ge\kappa_{1} d_s^2(d_n+1)$ it holds $q_{d_s}(t/N)\geq
d_s/2$ whenever $N-d_n-1\leq t\leq N-1$. Taking into account that
$\theta_N\to1$ as $N\to\infty$, it follows that for properly
selected $c_0(d_s,d_n)\geq\kappa_1d_s^2(d_n+1)$ and all $N\geq
c_0(d_s,d_n)$ it holds
%
\begin{eqnarray}
\label{eq2221} %
&& N-d_n-1\leq t\leq N-1\quad \Rightarrow
\quad | \bar{x}_t|=(-1)^t\bar{x}_t,
\nonumber
\\[-8pt]
\\[-8pt]
&&\min_{N-{d_n}-1\le t\le N-1} |\bar{x}_t|\ge|\bar{x}_{N-1}|/2
\ge \kappa_{2}\theta_N\lambda N^{-1/2}
{d_s}\geq\kappa_3\lambda N^{-1/2}d_s.
\nonumber
\end{eqnarray}
Beside this, for every ${d_s}$ and $N$, we clearly have $\bar{x}\in
\cS[\mathbf{w}]$ with $\mathbf{w}=\{\overbrace{\uppi,\ldots,\uppi
}^{{d_s}}\}$ (indeed, for all $t\in\bZ$, $((1+\Delta)^{d_s}\bar
{x})_t=(-1)^t\theta_N\lambda((1-\Delta)^{d_s} q)_t=0$).

2$^0$. Let $\overline{\mathbf{w}}=\{\overbrace
{0,\ldots,0}^{d_n}\}$, with $p_{\overline{\mathbf{w}}}(\zeta)=(1-\zeta)^{d_n}$.
Assuming $N\geq c_0(d_s,d_n)$, we have
\[
\bigl\llvert \bigl(p_{\overline{\mathbf{w}}}(\Delta)\bar{x} \bigr)_{N-1} \bigr
\rrvert \ge 2^{d_n} \min_{N-{d_n}-1\le t\le N-1} |\bar{x}_t|
\ge2^{d_n}\kappa _{3} \lambda N^{-1/2}
{d_s},
\]
so that for every $z\in\cS[\overline{\mathbf{w}}]$ we have
\[
\bigl| \bigl(p_{\overline{\mathbf{w}}}(\Delta)[{\bar{x}}-z] \bigr)_{N-1}\bigr|
\ge2^{d_n} \kappa_{3}\lambda N^{-1/2}
{d_s}.
\]
Since $N-1\geq d_n$, when taking into account that $|p_{\overline
{\mathbf{w}}}|_{1}=2^{d_n}$ we get for any $z\in\cS[\overline
{\mathbf{w}}]$
\[
2^{d_n}\bigl\|[\bar{x}-z]_0^{N-1}\bigr\|_\infty=|p_{\overline{\mathbf
{w}}}|_1
\bigl\|[\bar{x}-z]_0^{N-1}\bigr\|_\infty\geq\bigl|
\bigl(p_{\overline{\mathbf
{w}}}(\Delta)[\bar{x}-z] \bigr)_{N-1}\bigr|
\ge2^{d_n}\kappa_{3} \lambda N^{-1/2}
{d_s},
\]
whence
%
\begin{equation}
\label{whence12} \forall z\in\cS[\widehat{\mathbf{w}}]\dvtx \bigl \|[
\bar{x}-z]_0^{N-1} \bigr\| _\infty\geq
\kappa_{3} \lambda N^{-1/2} {d_s}.
\end{equation}
Now let us set $\lambda=\kappa_{4}\sqrt{\ln(1/\alpha)}$, with the
absolute constant $\kappa_4$ selected to ensure that $\lambda
<2{Q_{\cN}}(\alpha)$. The latter relation, due to $\lambda=\|
x_0^{N-1}\|_2$, ensures that
the hypotheses \emph{``observation (\ref{eqobs}) comes from $x\equiv
0$''} and \emph{``observation (\ref{eqobs}) comes from $x=\bar
{x}$''} cannot be distinguished $(1-\alpha)$-reliably. Thus, (\ref
{whence12}) implies the lower resolution bound $\kappa_3\kappa
_4d_sN^{-1/2}\sqrt{\ln(1/\alpha)}$.

\subsection{Proof of Proposition \texorpdfstring{\protect\ref{lowbndC}}{4.2}}\label{sec:proof41}
In the reasoning below, $c_i$ denote positive quantities depending
solely on $d=d_n$, and $\kappa_i$ denote positive absolute constants.
We start with proving the claim (ii).

1$^0$. Let us set
%
\begin{equation}
\label{f} \bar{f}_t=\sin \biggl({\uppi\over8}+
{\uppi\over4}{t\over N-1} \biggr), \qquad t\in\bZ.
\end{equation}
Assuming $c_0>40d$, let us fix an integer $\tau$ such that
%
\begin{equation}
\label{t} d-1 \leq\tau\leq(N-1)-20 d
\end{equation}
and set
%
\begin{equation}
\label{eq1234x} \gamma=\gamma_{\tau,N} = \epsilon\bar{f}^{-1}_\tau,
\qquad f=\gamma\bar{f}.
\end{equation}
By the definition of $f$ we have
%
\begin{eqnarray}
\label{sothat121a} 0\leq f_t\leq\epsilon\qquad \mbox{for } 0\leq t\leq
\tau,
\end{eqnarray}
and
%
\begin{eqnarray}
\label{sothat121b} \epsilon+\kappa_0N^{-1}(t-\tau)\epsilon\leq
f_t\leq\epsilon +\kappa_1N^{-1}(t-\tau)
\epsilon\qquad \mbox{for } \tau\leq t\leq N-1.
\end{eqnarray}
Let $p(\zeta)=(1-\zeta)^d=p_{\overline{\mathbf{w}}}(\zeta)$.
We clearly can find a sequence $x=\{x_t=a\cos({\uppi\over4}t+b)\}
_{t=-\infty}^\infty$ such that
$p(\Delta) x= f$, and, due to ${d_s}\geq2$, we have $x\in\cS_{{d_s}}$.
Further, let $\bar{z}\in\cS$ satisfy $\bar{z}_t=x_t$ for $0\leq
t\leq\tau$, and
\[
\bigl(p(\Delta)\bar{z} \bigr)_t = \lleft\{ %
\begin{array}
{l@{\qquad}l}f_t,&0\leq t\leq\tau,
\\
0,&t<0,
\\
\epsilon,&t\geq\tau.
\\
\end{array} %
\rright.
\]
Note that $\bar{z}$ is well defined due to $p(\Delta)x={f}$, and,
taking into account (\ref{eq1234x})--(\ref{sothat121b}), we
conclude that
$
\bar{z}\in\cS^\epsilon[\overline{\mathbf{w}}]$.

2$^0$. By the above construction, the sequence $\delta
=x-\bar{z}$ is such that $(p(\Delta)\delta)_t=0$ for $0\leq t\leq
\tau$, and $\kappa_0(t-\tau)N^{-1}\epsilon\leq(p(\Delta)\delta
)_t\leq\kappa_1(t-\tau)N^{-1}\epsilon$ when $\tau<t<N$, see (\ref
{sothat121b}). Besides, $\delta_t=0$ when $0\leq t\leq\tau$. By
evident reasons, these two observations combine with the first
inequality in (\ref{t}) to imply that
%
\begin{equation}
\label{implythat453} \bigl\|\delta_0^{N-1}\bigr\|_\infty\leq
c_4(N-\tau)^{d+1}N^{-1}\epsilon
\end{equation}
for some $c_4>0$ depending solely on $d$. We conclude that
\[
\bigl\|\delta_0^{N-1}\bigr\|_2\leq\bigl\|\delta_0^{N-1}
\bigr\|_\infty\sqrt{N-\tau }\leq c_4(N-\tau)^{d+3/2}N^{-1}
\epsilon.
\]
Now note that the {hypotheses ``observation (\ref{eqobs}) comes from
$x= \bar{z}$'' and ``observation (\ref{eqobs}) comes from $x=\bar
{z}+\delta$'' cannot be distinguished $(1-\alpha)$-reliably unless $\|
\delta_0^{N-1}\|_2\ge2{Q_{\cN}}(\alpha)> \sqrt{\kappa_3\ln(1/
\alpha)}$}.
Equipped with this $\kappa_3$ and with $c_4$ participating in (\ref
{implythat453}), let us set
%
\begin{equation}
\label{nu} \nu=\Floor \bigl( \bigl(\kappa_3\ln(1/\alpha
)N^2c_4^{-2}\epsilon^{-2}
\bigr)^{{1/(2d+3)}} \bigr).
\end{equation}
It is immediately seen that with properly chosen positive $c_1,c_3$
depending solely on $d$, and with $\epsilon$ satisfying (\ref{range}),
we have $20d<\nu<N-d$, so that setting
\[
\tau=N-\nu,
\]
we ensure (\ref{t}). From now on, we assume that $c_1,c_3$ are as
needed in the latter conclusion. With the just defined $\tau$, we have
$
\|\delta_0^{N-1}\|_2^2\leq\kappa_3\ln(1/\alpha)$,
meaning that $x$ and $\bar{z}$ cannot be distinguished $(1-\alpha)$-reliably.

3$^0$. To prove the claim (ii) it now suffices to show that a
properly chosen $c_5>0$, depending solely on $d$,
%
\begin{equation}
\label{deviation} \forall z\in\cS^\epsilon[\overline{\mathbf{w}}]\dvtx
\bigl\|[x-z]_0^{N-1}\bigr\| _\infty\geq c_5(N-
\tau)^{d+1}N^{-1}\epsilon=c_5\nu
^{d+1}N^{-1}\epsilon.
\end{equation}
Indeed, given $z\in\cS^\epsilon$, we have
\[
\tau\leq t<N\quad \Rightarrow\quad \bigl(p(\Delta)[x-z] \bigr)_t=
\bigl(p(\Delta)x \bigr)_t- \bigl(p(\Delta )z \bigr)_t\geq
f_t- \epsilon\geq\kappa_0N^{-1}\epsilon(t-
\tau),
\]
with concluding inequality given by (\ref{sothat121b}). Setting
$\theta
=\lfloor(N-1-\tau)/2\rfloor$, the sequence
\[
s_t=(x-z)_{t-\tau},\qquad t\in\bZ
\]
satisfies
%
\begin{equation}
\label{satisfies22} \bigl(p(\Delta)s \bigr)_t\geq\kappa_2
\theta N^{-1}\epsilon, \qquad \theta\leq t\leq 2\theta,
\end{equation}
and, by the right inequality in (\ref{t}), $\theta\geq10d$. Setting
$k=\lfloor(\theta-d)/d\rfloor\geq2$, consider the polynomial
\[
q(\zeta)= \bigl(1-\zeta^k \bigr)^d=(1-
\zeta)^d\underbrace{ \bigl(1+\zeta+\cdots+\zeta ^{k-1}
\bigr)^d}_{v(\zeta)}=(1-\zeta)^d\sum
_{j=0}^{d(k-1)}v_j\zeta^j;
\]
where, clearly,
%
\begin{eqnarray}
v_j\geq0,\qquad \sum_jv_j=k^d.
\label{clearly}
\end{eqnarray}
Let now $r=q(\Delta)s$. Taking into account that $|q(\cdot)|_{1}=
2^d$, we have
%
\begin{equation}
\label{wehave3334} r_{2\theta}\leq2^d\bigl\|s_\theta^{2\theta}
\bigr\|_\infty\le2^d\bigl\| [x-z]_0^{N-1}
\bigr\|_\infty,
\end{equation}
since, by construction,
$\|[x-z]_0^{N-1}\|_\infty\geq\|s_\theta^{2\theta}\|_\infty$.
On the other hand, $r=v(\Delta)u$, $u=p(\Delta)s$, so that
\[
r_{2\theta}=\sum_{j=0}^{d(k-1)}v_ju_{2\theta-j}.
\]
By (\ref{satisfies22}), we have $u_{2\theta-j}\geq\kappa_2\theta
N^{-1}\epsilon$ for $0\leq j\leq d(k-1)$ (note that $d(k-1)<\theta$),
and by (\ref{clearly}),
\[
r_{2\theta}\geq\kappa_2\theta N^{-1}\epsilon\sum
_j{v_j}=k^d\kappa
_2\theta N^{-1}\epsilon.
\]
Combining the latter inequality with (\ref{wehave3334}) we come to $\|
[\bar{x}-z]_0^{N-1}\|_\infty\geq\kappa_22^{-d}k^d\theta
N^{-1}\epsilon$. Recalling that by construction $k>\kappa_{10}(N-\tau
)/d$, $\theta\geq\kappa_{11}(N-\tau)$, we arrive at (\ref
{deviation}). (ii) is proved.

4$^0$. It remains to prove (i).
Note that when $\epsilon<c_{12}N^{-1/2}\sqrt{\ln(1/\alpha)}$, the
conclusion in (i) is readily given by a straightforward modification of
the reasoning in Section~\ref{sec:proof40}. Note that for the time
being the only restriction on the lower bound $c_0$ on $N$, see (\ref
{proplet}), was that $c_0\geq40d$, see the beginning of item 1$^0$.
Now let us also assume that $N$ is large enough to ensure that
$c_{12}N^{-1/2}\geq c_3N^{-d-1/2}$ (which still allows to choose $c_0$
as a function of $d=d_n$ only). With the resulting $c_0$, the range of
values of $\epsilon$ for which we have justified the conclusion in (i)
covers the corresponding range of $\epsilon$ allowed by the premise of (i).

\subsection{Proof of Proposition \texorpdfstring{\protect\ref{lowbndD}}{4.3}}\label{sec:proof42}
The proof of the statement (i) is given by a straightforward
modification of the reasoning in Section~\ref{sec:proof40}. Let us
prove (ii).
Let $\mathbf{w}=\{\overbrace{\uppi,\ldots,\uppi}^{d_s}\}$, and $0\le\tau
< N-d_n-1$, so that $\cS^{N,\epsilon_s}[{\mathbf{w}}]$ contains the sequence
$\bar{x}=\{\bar{x}_t=c_1\epsilon_s(-1)^t (t-\tau)_+^{d_s}\}
_{t=-\infty}^\infty$ with $0< c_1=c_1(d_s)$ small enough. Then for
$\overline{\mathbf{w}}=\{0,\ldots,0\}$ ($d_n$ zeros), with $p_{\overline
{\mathbf{w}}}(\zeta)=(1-\zeta)^{d_n}$, we have for any $z\in\cS
^{N,\epsilon_n}[\overline{\mathbf{w}}]$:
\begin{eqnarray*}
\bigl| \bigl(p_{\overline{\mathbf{w}}}(\Delta)[{\bar{x}}-z] \bigr)_{N-1} \bigr|&\ge& \bigl|
\bigl(p_{\overline{\mathbf{w}}}(\Delta)\bar{x} \bigr)_{N-1}\bigr |- \bigl|
\bigl(p_{\overline
{\mathbf{w}}}(\Delta)z \bigr)_{N-1}\bigr|\ge 2^{d_n}|
\bar{x}_{N-d-1} |-\epsilon_n
\\
&\geq& c_2 \epsilon_s(N-\tau)^{d_s}-
\epsilon_n,
\end{eqnarray*}
where $c_2>0$ depends only on $d_s$ and $d_n$. Therefore, when denoting
$\nu=N-\tau$, for every $z\in\cS^{N,\epsilon_n}[\overline{\mathbf
{w}}]$ we have
\[
2^{d_n}\bigl\|[\bar{x}-z]_0^{N-1}\bigr\|_\infty=|p_{\overline{\mathbf
{w}}}|_1
\bigl\|[\bar{x}-z]_0^{N-1}\bigr\|_\infty\geq\bigl|
\bigl(p_{\overline{\mathbf
{w}}}(\Delta)[\bar{x}-z] \bigr)_{N-1}\bigr|\ge
c_2 \epsilon_s\nu ^{d_s}-\epsilon_n,
\]
whence,
%
\begin{equation}
\label{whence121} \forall z\in\cS^{N,\epsilon_n}[\overline{\mathbf{w}}]\dvtx \bigl\|[
\bar {x}-z]_0^{N-1}\bigr\|_\infty\geq2^{-d_n}
\bigl(c_2\epsilon_s\nu^{d_s}-
\epsilon_n \bigr).
\end{equation}
It is immediately seen that with properly selected positive
$c_i=c_i(d_n,d_s)$, $i=3,4,5,6$, assuming
%
\begin{equation}
\label{assuming776} \epsilon_s N^{d_s+\sfrac{1}{2}}\geq c_3
\sqrt{ \ln(1/\alpha)} \quad \mbox{and}\quad \epsilon _s\leq
c_4\sqrt{\ln(1/ \alpha)}
\end{equation}
and selecting $\tau$ according to
\[
\nu=N-\tau=\Floor \bigl(c_5 \bigl[\epsilon_s^{-2}
\ln (1/\alpha) \bigr]^{{1/(2d_s+1)}} \bigr),
\]
we ensure that $0\leq\tau<N-d_n-1$ and
\begin{eqnarray*}
&& \mbox{(a)}\quad \tau<N-d_n-1,
\\
&&\mbox{(b)}\quad \bigl\| \bar{x}_{0}^{N-1}\bigr\|_2<
2{Q_{\cN}}( \alpha),
\\
&&\mbox{(c)}\quad 2^{-d_n} \bigl(c_2\epsilon_s
\nu^{d_s}-\epsilon_n \bigr)\geq\bar{\rho }:=c_6
\epsilon_s^{1/(2d_s+1)} \bigl(\ln(1/\alpha) \bigr)^{{d_s/(2d_s+1)}}
\end{eqnarray*}
(when verifying (c), take into account that $0\leq\epsilon_n\leq
\epsilon_s$). By (b), the hypotheses ``observation (\ref{eqobs})
comes from $x\equiv0$'' and ``observation (\ref{eqobs}) comes from
$x=\bar{x}$'' cannot be distinguished $(1-\alpha)$-reliably, while
(\ref{whence121}), the already established inclusion $\bar{x}\in\cS
^{N,\epsilon_s}[{\mathbf{w}}]$ and (c) imply that $\bar{x}$ obeys
the hypothesis $H_1(\bar{\rho})$ associated with the problem $(N_2)$.
The bottom line is that in the case of (\ref{assuming776}), $\bar
{\rho}$ defined in (c) is a lower bound on $\rho_*(\alpha)$.

\begin{appendix}

\section*{Appendix: Proof of Lemma \texorpdfstring{\protect\ref{MainLemma}}{6.1}}\label{app}
\renewcommand{\theequation}{\arabic{equation}}
\setcounter{equation}{43}
1$^0$. Let $\lambda\in\C$, $0<|\lambda|\leq1$, let
$\epsilon\in(0,1)$, and let $n\geq1$ be an integer. Setting $\delta
=1-\epsilon$, consider the polynomials
%
\begin{eqnarray}
\label{polynomials} f_k( \zeta)&=&f_k({\zeta};\lambda,
\epsilon):=(1-\lambda\zeta)\sum_{\ell=0}^k(
\delta\lambda\zeta)^\ell= {1-\lambda{\zeta}\over
1-\delta\lambda{\zeta}} \bigl[1-[\delta
\lambda\zeta]^{k+1} \bigr],\qquad k=0,1,\ldots
\nonumber
\\[-8pt]
\\[-8pt]
p_n(\zeta)&=&p_n({\zeta};\lambda,\epsilon):=
{1\over n}\sum_{k=0}^{n-1}f_k(
\zeta)=(1-\lambda\zeta)\sum_{\ell=0}^{n-1}[
\delta \lambda\zeta]^\ell[1-\ell/n].
\nonumber
\end{eqnarray}
Observe that
\[
f_k(0)=1,\qquad f_k(1/\lambda)=0,
\]
whence
%
\begin{equation}
\label{eq3} p_n(0)=1,\qquad p_n(1/\lambda)=0.
\end{equation}
%

1.1$^0$. Let us bound from above the uniform norm $\|
p_n(\cdot)\|_\infty$ of $p_n$ on the unit circle. We have $p_n(\zeta
)=r_n(\lambda\zeta)$, where
\begin{eqnarray*}
r_n(\zeta)&=&{1\over n}(1- \zeta)\sum
_{k=0}^{n-1}{1-(\delta\zeta
)^{k+1}\over1-\delta\zeta}=(1- \zeta)
{n-\delta\zeta{((1-[\delta
\zeta]^{n})/(1-\delta\zeta))}\over n(1-\delta\zeta)}
\\
&=&\underbrace{{1-\zeta\over1-\delta\zeta}}_{g_n(\zeta
)}\underbrace{
{n(1-\delta\zeta)-\delta\zeta+[\delta\zeta
]^{n+1}\over n(1-\delta\zeta)}}_{h_n(\zeta)}.
\end{eqnarray*} %
Since $|\lambda|\leq1$, we have $\|p_n(\cdot)\|_\infty\leq\|
g_n(\cdot)\|_\infty\|h_n(\cdot)\|_\infty$. When $|\zeta|=1$, we have
%
\begin{equation}
\label{conc1} \bigl|g_n(\zeta)\bigr|={1\over\delta}
{|1-\zeta|\over|\delta^{-1}-
\zeta|}\leq{1\over\delta}.
\end{equation}
If we set $
\zeta=\cos(\phi)+\imath\sin(\phi)$, and
$\epsilon=\theta/ n$ and $\delta=1-\epsilon=1-{\theta/ n}$
with some $\theta\in(0,n)$,
we obtain
\begin{eqnarray*}
\bigl|h_n(\zeta)\bigr|&\leq&{|n(1-\delta\zeta)-\delta\zeta|+\delta
^{n+1}\over n|1-\delta\lambda|}
\nonumber
\\
&=&{\sqrt{[n-(n+1)\delta\cos(\phi)]^2+(n+1)^2\delta^2\sin^2(\phi
)}+\delta^{n+1}\over
n\sqrt{[1-\delta\cos(\phi)]^2+\delta^2\sin^2(\phi)}}
\nonumber
\\
&=&{\sqrt{n^2+(n+1)^2\delta^2\cos^2(\phi)-2n(n+1)\delta\cos(\phi
)+(n+1)^2\delta^2\sin^2(\phi)}+\delta^{n+1}\over
n\sqrt{1+\delta^2\cos^2(\phi)-2\delta\cos(\phi)+\delta^2\sin
^2(\phi)}}
\nonumber
\\
&=&{\sqrt{n^2+(n+1)^2-2n(n+1)\delta\cos(\phi)}+\delta^{n+1}\over
n\sqrt{1+\delta^2-2\delta\cos(\phi)}}
\nonumber
\\
&=&{\sqrt{[n-(n+1)\delta)]^2+2n(n+1)\delta[1-\cos(\phi)]}+\delta
^{n+1}\over n\sqrt{[1-\delta]^2+2\delta[1-\cos(\phi)]}}
\nonumber
\\
&=&{\sqrt{ [\theta{((n+1)/ n)}-1 ]^2+2(n+1)(n-\theta
)[1-\cos(\phi)]}+\delta^{n+1}\over n\sqrt{{\theta^2/
n^2}+2n^{-1}(n-\theta)[1-\cos(\phi)]}}
\nonumber
\\
&=&{\sqrt{ [\theta{((n+1)/ n)}-1 ]^2+2(n+1)(n-\theta
)[1-\cos(\phi)]}+\delta^{n+1}\over\sqrt{\theta^2+2n(n-\theta
)[1-\cos(\phi)]}} ={\sqrt{\alpha+\beta t} +\gamma\over\sqrt{\mu+\nu t}},
\end{eqnarray*}
where
%
\begin{eqnarray}
\label{alphaetal} %
\alpha&=& \bigl[\theta{\bigl((n+1)/ n\bigr)}-1
\bigr]^2,\qquad \beta=2(n+1) (n-\theta), \qquad \gamma=
\delta^{n+1},
\nonumber
\\[-8pt]
\\[-8pt]
\mu&=& \theta^2,\qquad \nu=2n(n-\theta), \qquad t=1-\cos(\phi).
\nonumber
\end{eqnarray}
It is immediately seen that with positive $\alpha,\beta,\gamma,\mu$
and $\nu$ such that
%
\begin{eqnarray}
\mu/\nu> \alpha/\beta, \label{(!)123}
\end{eqnarray}
the maximum of the function ${\sqrt{\alpha+\beta t}+\gamma\over
\sqrt{\mu+\nu t}}$ over $t$ such that $\alpha+\beta t\geq0$ is achieved
when $\alpha+\beta t=\sqrt{{\beta\mu-\alpha\nu\over\gamma\nu
}}$ and is equal to $\sqrt{{\beta\over\nu}}\sqrt{1+{\gamma^2\nu
\over\beta\mu-\alpha\nu}}$. Now, assume that
%
\begin{equation}
\label{assumingthat} 1<\theta< n.
\end{equation}
Then ${n+1\over n}\theta^2> [\theta{n+1\over n}-1]^2$,
so that
the parameters $\alpha,\ldots,\nu$ defined in (\ref{alphaetal}) satisfy
(\ref{(!)123}).
We conclude that
%
\begin{equation}
\label{conclude} \bigl\|h_n(\cdot)\bigr\|_\infty\leq\sqrt{
{n+1\over n}}\sqrt{1+{[1-\theta
/n]^{2(n+1)}\over{{((n+1)/ n)}\theta^2- [\theta{((n+1)/
n)}-1 ]^2}}}.
\end{equation}
Note that
\[
{n+1\over n}\theta^2- \biggl[\theta{n+1\over n}-1
\biggr]^2=2\theta{n+1\over
n}-\theta^2
{n+1\over n^2}-1\ge1
\]
when $\theta$ satisfies
%
\begin{equation}
\label{assumingmore} {2n\over n+1}\leq\theta<n,
\end{equation}
and (\ref{conclude}) implies in this case that
\[
\bigl\|h_n(\cdot)\bigr\|_\infty\leq\sqrt{{n+1\over n}}
\sqrt{1+[1-\theta /n]^{2(n+1)}}\leq\exp \biggl\{{1\over2n}+
{1\over2}\e^{-2\theta
} \biggr\}.
\]
The latter bound combines with (\ref{conc1}) to imply that in the case
of (\ref{assumingmore}) (recall that $1-\delta=\epsilon=\theta/n$)
we get for all $|\lambda|\le1$:
%
\begin{eqnarray}
\label{boundinf} \max_{|z|\leq1}\bigl|p_n(z;\lambda,
\epsilon)\bigr|&\leq& {1\over1-\theta
/n}\exp \biggl\{{1\over2n}+
{1\over2}\e^{-2\theta} \biggr\}
\nonumber
\\[-8pt]
\\[-8pt]
&\leq&\exp \biggl\{{\theta\over n-\theta}+{1\over2n}+
{1\over2}\e ^{-2\theta} \biggr\}.
\nonumber
\end{eqnarray}
%

1.2$^0$.
Now let us bound from above $\|1-p_n(\cdot)\|_2$.
We have
\begin{eqnarray*}
p_n(\zeta)&=&1+\sum_{\ell=1}^{n-1}
\biggl(\delta^\ell \biggl[1- {\ell
\over n} \biggr]-
\delta^{\ell-1} \biggl[1- {\ell-1\over n} \biggr] \biggr)[\lambda
\zeta]^\ell- {1\over n}\delta^{n-1}[\lambda\zeta
]^n
\\
&=&1+\sum_{\ell=1}^{n-1}\delta^{\ell-1}
\biggl([1-\epsilon] \biggl[1-{\ell\over n} \biggr]- \biggl[1-
{\ell-1\over n} \biggr] \biggr)[\lambda\zeta]^\ell-
{1\over n}\delta^{n-1}[\lambda\zeta]^n
\\
&=&1- \Biggl(\sum_{\ell=1}^{n-1}
\delta^{\ell-1} \biggl[\epsilon \biggl[1-{\ell\over n} \biggr]+
{1\over n} \biggr] [\lambda\zeta]^\ell+
{1\over
n}\delta^{n-1}[\lambda\zeta]^n \Biggr).
\end{eqnarray*} %
Taking into account that $|\lambda|\leq1$, we conclude that
%
\begin{eqnarray}
\label{eq2} \bigl\|1-p_n(\cdot)\bigr\|_2&\leq& \Biggl\|\sum
_{\ell=1}^{n-1}(1-\epsilon)^{\ell
-1}\epsilon
\biggl[1-{\ell\over n} \biggr] [\lambda\zeta]^\ell
\Biggr\|_2+{1\over n}\Biggl\| \sum_{\ell=1}^{n}(1-
\epsilon)^{\ell-1}[\lambda\zeta]^\ell\Biggr\| _2
\nonumber
\\
&\leq&\sqrt{\sum_{\ell=1}^n
\epsilon^2(1-\epsilon)^{2(\ell
-1)}}+\sqrt{1/n}
\\
&\leq&\sqrt{{\epsilon^2\over1-(1-\epsilon)^2}}+\sqrt{1/n}\leq \sqrt{\epsilon}+\sqrt
{1\over n}\leq2\sqrt{\max \biggl[\epsilon ,{1\over n}
\biggr]}.
\nonumber
\end{eqnarray}
%

2$^0$. Let $n=\lfloor m/d\rfloor$ and $\lambda_\ell
=\exp\{-\imath\upsilon_\ell\}$, $1\leq\ell\leq d$, where $m$, $d$ and
$\upsilon_1,\ldots,\upsilon_d$ are as described in the premise of the
lemma. Let us set
\begin{eqnarray*}
\theta&=&\max \biggl[2,{1\over2}\ln(2d) \biggr]\leq3\ln(2d),
\qquad \epsilon={\theta
\over n},
\\
q(\zeta)&=&1-p_n({\zeta};\lambda_1,\epsilon)\cdot
p_n({\zeta };\lambda_2,\epsilon)\cdot\ldots \cdot
p_n({\zeta;}\lambda_d,\epsilon).
\end{eqnarray*}
%

2.1$^0$. Observe that by (\ref{mislargeenough}) we have
%
\begin{equation}
\label{nislarge} n\geq n(d):=\Ceil(5d\theta)=\Ceil \biggl(5d\max \biggl[2,
{1\over
2}\ln(2d) \biggr] \biggr).
\end{equation}
Our choice of $\theta$ and $n$ ensures (\ref{assumingmore}), so that
$0<\epsilon<1$, and by (\ref{boundinf}) we also have
%
\begin{equation}
\label{eq112} \max_{|z|\leq1}\bigl|p_n(z;
\lambda_\ell,\epsilon)\bigr|\leq\exp \biggl\{ {\theta\over n-\theta}+
{1\over2n}+{1\over2}\mathrm{e}^{-2\theta} \biggr
\} \leq \e^{1/d},
\end{equation}
where the concluding inequality is readily given by the choice of
$\theta$ and by (\ref{nislarge}). Recall that by (\ref
{polynomials}), $p_n({\zeta};\lambda_\ell,\epsilon)$ is divisible
by $(1-\lambda_\ell\zeta)$; when setting\vspace*{1pt}
\[
r_\ell(\zeta)={p_n({\zeta};\lambda_\ell,\epsilon)\over1-\lambda
_\ell\zeta}=\sum
_{\ell=0}^{n-1} \bigl[(1-\epsilon)\lambda_\ell
\zeta \bigr]^\ell[1-\ell/n],
\]
and $r(\zeta)=\prod_{\ell=1}^d r_\ell(\zeta)$, we clearly have
$r(0)=1$, $\deg r\leq d(n-1)\leq m-d$, and\vspace*{1pt}
\[
1-q(\zeta)=p_{{\mathbf{u}}}(\zeta)\prod_{\ell=1}^d
r_\ell(\zeta )=p_{{\mathbf{u}}}(\zeta)r(\zeta),
\]
as required in
(\ref{divisible}) (note that $q$ is a real polynomial due to
${{\mathbf
{u}}}\in\Omega_d$).

2.2$^0$. By (\ref{eq112}) we have $\|p_n(\cdot;\epsilon
,\lambda_k)\|_\infty\leq\e^{1/d}$, while (\ref{eq2}) says that\vspace*{1pt}
\[
\bigl\|p_n(\cdot;\lambda_k,\epsilon)-1\bigr\|_2\leq3
\sqrt{\ln(2d)/n}.
\]
We have\vspace*{1pt}
\[
1-\prod_{k=1}^{\ell+1}p_n(\cdot;
\lambda_k,\epsilon)= \Biggl[1-\prod_{k=1}^{\ell}p_n(
\cdot;\lambda_k,\epsilon) \Biggr]p_n(\cdot ;
\lambda_{\ell+1},\epsilon)+ \bigl[1-p_n(\cdot;
\lambda_{\ell+1},\epsilon) \bigr].
\]
When denoting $\alpha_\ell=\|1-\prod_{k=1}^{\ell}p_n(\cdot;\lambda
_k,\epsilon)\|_2$ for $\ell=1,2,\ldots,d$ and setting $\alpha_0=0$, we get\vspace*{1pt}
\[
\alpha_{\ell+1}\leq\alpha_\ell\bigl\|p_n(\cdot;
\lambda_{\ell
+1},\epsilon)\bigr\|_\infty+\bigl\|p_n(\cdot;
\lambda_{\ell+1},\epsilon)-1\bigr\| _2\leq\alpha_\ell
\e^{1/d}+3\sqrt{\ln(2d)/n}, \qquad 0\leq\ell<d.
\]
It follows that $\alpha_\ell\leq3\ell\sqrt{\ln(2d)\over n}\e
^{\ell/d}$ when $\ell\leq d$, which implies
\[
\bigl\|q(\cdot)\bigr\|_2\leq3\mathrm{e} d\sqrt{\ln(2d)/n}\leq3\e
d^{3/2}\sqrt {\ln
(2d)\over m},
\]
and (\ref{propofq}) follows.

\end{appendix}




\printhistory


\begin{thebibliography}{20}


\bibitem{Chiu}
\begin{barticle}[mr]
\bauthor{\bsnm{Chiu},~\bfnm{Shean-Tsong}\binits{S.-T.}}
(\byear{1989}).
\btitle{Detecting periodic components in a white {G}aussian time series}.
\bjournal{J. Roy. Statist. Soc. Ser. B}
\bvolume{51}
\bpages{249--259}.
\bid{issn={0035-9246}, mr={1007457}}
\end{barticle}
\bptok{imsref}%
\endbibitem

\bibitem{Davies87}
\begin{barticle}[mr]
\bauthor{\bsnm{Davies},~\bfnm{Robert~B.}\binits{R.B.}}
(\byear{1987}).
\btitle{Hypothesis testing when a nuisance parameter is present only under the alternative}.
\bjournal{Biometrika}
\bvolume{74}
\bpages{33--43}.
\bid{issn={0006-3444}, mr={0885917}}
\end{barticle}
\bptok{imsref}%
\endbibitem

\bibitem{Dju96}
\begin{barticle}[auto:STB|2014/02/12|14:17:21]
\bauthor{\bsnm{Djuric},~\bfnm{P.~M.}\binits{P.M.}}
(\byear{1996}).
\btitle{A model selection rule for sinusoids in white Gaussian noise}.
\bjournal{IEEE Trans. Signal Process.}
\bvolume{44}
\bpages{1744--1751}.
\end{barticle}
\bptok{imsref}%
\endbibitem

\bibitem{Fish29}
\begin{barticle}[auto:STB|2014/02/12|14:17:21]
\bauthor{\bsnm{Fisher},~\bfnm{R.~A.}\binits{R.A.}}
(\byear{1929}).
\btitle{Test of significance in harmonic analysis}.
\bjournal{Proc. R. Soc. Lond. Ser. A}
\bvolume{125}
\bpages{54--59}.
\end{barticle}
\bptok{imsref}%
\endbibitem

\bibitem{GoNem97}
\begin{barticle}[mr]
\bauthor{\bsnm{Goldenshluger},~\bfnm{Alexander}\binits{A.}} \AND
\bauthor{\bsnm{Nemirovski},~\bfnm{Arkadi}\binits{A.}}
(\byear{1997}).
\btitle{Adaptive de-noising of signals satisfying differential inequalities}.
\bjournal{IEEE Trans. Inform. Theory}
\bvolume{43}
\bpages{872--889}.
\bid{doi={10.1109/18.568698}, issn={0018-9448}, mr={1454226}}
\end{barticle}
\bptok{imsref}%
\endbibitem

\bibitem{Han70}
\begin{bbook}[mr]
\bauthor{\bsnm{Hannan},~\bfnm{E.~J.}\binits{E.J.}}
(\byear{1970}).
\btitle{Multiple Time Series}.
\blocation{New York}:
\bpublisher{Wiley}.
\bid{mr={0279952}}
\end{bbook}
\bptok{imsref}%
\endbibitem

\bibitem{Han93}
\begin{bincollection}[mr]
\bauthor{\bsnm{Hannan},~\bfnm{E.~J.}\binits{E.J.}}
(\byear{1993}).
\btitle{Determining the number of jumps in a spectrum}.
In \bbooktitle{Developments in Time Series Analysis}
(\beditor{\bfnm{T. Subba}\binits{T.S.} \bsnm{Rao}}, ed.)
\bpages{127--138}.
\blocation{London}:
\bpublisher{Chapman \& Hall}.
\bid{mr={1292263}}
\end{bincollection}
\bptok{imsref}%
\endbibitem

\bibitem{JuNem1}
\begin{barticle}[mr]
\bauthor{\bsnm{Juditsky},~\bfnm{Anatoli}\binits{A.}} \AND
\bauthor{\bsnm{Nemirovski},~\bfnm{Arkadi}\binits{A.}}
(\byear{2009}).
\btitle{Nonparametric denoising of signals with unknown local structure, {I}: {O}racle inequalities}.
\bjournal{Appl. Comput. Harmon. Anal.}
\bvolume{27}
\bpages{157--179}.
\bid{doi={10.1016/j.acha.2009.02.001}, issn={1063-5203}, mr={2543191}}
\end{barticle}
\bptok{imsref}%
\endbibitem

\bibitem{JuNem2}
\begin{barticle}[mr]
\bauthor{\bsnm{Juditsky},~\bfnm{Anatoli}\binits{A.}} \AND
\bauthor{\bsnm{Nemirovski},~\bfnm{Arkadi}\binits{A.}}
(\byear{2010}).
\btitle{Nonparametric denoising signals of unknown local structure, {II}: Nonparametric function recovery}.
\bjournal{Appl. Comput. Harmon. Anal.}
\bvolume{29}
\bpages{354--367}.
\bid{doi={10.1016/j.acha.2010.01.003}, issn={1063-5203}, mr={2672231}}
\end{barticle}
\bptok{imsref}%
\endbibitem

\bibitem{KH94}
\begin{barticle}[mr]
\bauthor{\bsnm{Kavalieris},~\bfnm{L.}\binits{L.}} \AND
\bauthor{\bsnm{Hannan},~\bfnm{E.~J.}\binits{E.J.}}
(\byear{1994}).
\btitle{Determining the number of terms in a trigonometric regression}.
\bjournal{J. Time Series Anal.}
\bvolume{15}
\bpages{613--625}.
\bid{doi={10.1111/j.1467-9892.1994.tb00216.x}, issn={0143-9782}, mr={1312325}}
\end{barticle}
\bptok{imsref}%
\endbibitem

\bibitem{NK11}
\begin{barticle}[mr]
\bauthor{\bsnm{Nadler},~\bfnm{Boaz}\binits{B.}} \AND
\bauthor{\bsnm{Kontorovich},~\bfnm{Aryeh}\binits{A.}}
(\byear{2011}).
\btitle{Model selection for sinusoids in noise: Statistical analysis and a new penalty term}.
\bjournal{IEEE Trans. Signal Process.}
\bvolume{59}
\bpages{1333--1345}.
\bid{doi={10.1109/TSP.2011.2105482}, issn={1053-587X}, mr={2807736}}
\end{barticle}
\bptok{imsref}%
\endbibitem

\bibitem{Nem81}
\begin{barticle}[mr]
\bauthor{\bsnm{Nemirovski{\u\i}},~\bfnm{A.~S.}\binits{A.S.}}
(\byear{1981}).
\btitle{Prediction under conditions of indeterminacy}.
\bjournal{Problems Inform. Transmission}
\bvolume{17}
\bpages{73--83}
(\bnote{in Russian}).
\bid{mr={0673748}}
\end{barticle}
\bptok{imsref}%
\endbibitem

\bibitem{Nem92}
\begin{bincollection}[mr]
\bauthor{\bsnm{Nemirovski{\u\i}},~\bfnm{A.~S.}\binits{A.S.}}
(\byear{1992}).
\btitle{On nonparametric estimation of functions satisfying differential inequalities}.
In \bbooktitle{Topics in Nonparametric Estimation}
(\beditor{\bfnm{R.}\binits{R.} \bsnm{Khasminskii}}, ed.).
\bseries{Adv. Soviet Math.}
\bvolume{12}
\bpages{7--43}.
\blocation{Providence, RI}:
\bpublisher{Amer. Math. Soc.}
\bid{mr={1191690}}
\end{bincollection}
\bptok{imsref}%
\endbibitem

\bibitem{Pis73}
\begin{barticle}[auto:STB|2014/02/12|14:17:21]
\bauthor{\bsnm{Pisarenko},~\bfnm{V.~F.}\binits{V.F.}}
(\byear{1973}).
\btitle{The retrieval of harmonics from a covariance function}.
\bjournal{Geophys. J. Roy. Astron. Soc.}
\bvolume{33}
\bpages{347--366}.
\end{barticle}
\bptok{imsref}%
\endbibitem

\bibitem{HanQui}
\begin{bbook}[mr]
\bauthor{\bsnm{Quinn},~\bfnm{B.~G.}\binits{B.G.}} \AND
\bauthor{\bsnm{Hannan},~\bfnm{E.~J.}\binits{E.J.}}
(\byear{2001}).
\btitle{The Estimation and Tracking of Frequency}.
\bseries{Cambridge Series in Statistical and Probabilistic Mathematics}
\bvolume{9}.
\blocation{Cambridge}:
\bpublisher{Cambridge Univ. Press}.
\bid{mr={1813156}}
\end{bbook}
\bptok{imsref}%
\endbibitem

\bibitem{QK94}
\begin{barticle}[auto:STB|2014/02/12|14:17:21]
\bauthor{\bsnm{Quinn},~\bfnm{B.~G.}\binits{B.G.}} \AND
\bauthor{\bsnm{Kotsookos},~\bfnm{P.~J.}\binits{P.J.}}
(\byear{1994}).
\btitle{Threshold behavior of the maximum likelihood estimator of frequency}.
\bjournal{IEEE Trans. Signal Process.}
\bvolume{42}
\bpages{3291--3294}.
\end{barticle}
\bptok{imsref}%
\endbibitem

\bibitem{RB74}
\begin{barticle}[auto:STB|2014/02/12|14:17:21]
\bauthor{\bsnm{Rife},~\bfnm{D.~C.}\binits{D.C.}} \AND
\bauthor{\bsnm{Boorstyne},~\bfnm{R.~R.}\binits{R.R.}}
(\byear{1974}).
\btitle{Single-tone parameter estimation from discrete-time observations}.
\bjournal{IEEE Trans. Inform. Theory}
\bvolume{20}
\bpages{591--598}.
\end{barticle}
\bptok{imsref}%
\endbibitem

\bibitem{mus1}
\begin{barticle}[auto:STB|2014/02/12|14:17:21]
\bauthor{\bsnm{Schmidt},~\bfnm{R.~O.}\binits{R.O.}}
(\byear{1986}).
\btitle{Multiple emitter location and signal parameter estimation}.
\bjournal{IEEE Trans. Antennas and Propagation}
\bvolume{34}
\bpages{276--280}.
\end{barticle}
\bptok{imsref}%
\endbibitem


\bibitem{mus2}
\begin{bbook}[auto:STB|2014/02/12|14:17:21]
\bauthor{\bsnm{Stoica},~\bfnm{P.}\binits{P.}} \AND
\bauthor{\bsnm{Moses},~\bfnm{R.~L.}\binits{R.L.}}
(\byear{1997}).
\btitle{Introduction to Spectral Analysis}.
\blocation{New York}:
\bpublisher{Prentice-Hall}.
\end{bbook}
\bptok{imsref}%
\endbibitem

\bibitem{Whittle}
\begin{barticle}[mr]
\bauthor{\bsnm{Whittle},~\bfnm{P.}\binits{P.}}
(\byear{1954}).
\btitle{The statistical analysis of a seiche record}.
\bjournal{J. Marine Res.}
\bvolume{13}
\bpages{76--100}.
\bid{mr={0078600}}
\end{barticle}
\bptok{imsref}%
\endbibitem

\end{thebibliography}
\end{document}